\documentclass[11pt]{article}

\usepackage{graphicx}
\usepackage{amssymb}
\usepackage{amscd}
\usepackage{amsmath}
\usepackage{enumerate}

\newtheorem{df}{Definition}[section]
\newtheorem{thm}{Theorem}
\newtheorem{pp}[df]{Proposition}
\newtheorem{lemma}[df]{Lemma}

\newtheorem{corollary}[df]{Corollary}

\newtheorem{prf}{{\it Proof.}}

\usepackage{picins}

\makeatletter

\@addtoreset{equation}{section}

\@addtoreset{figure}{section}
\makeatother

\newcommand{\qed}{\hbox{\rule[0pt]{3pt}{6pt}}}

\title{{\bf\large Rational Maps and Maximum Likelihood Decodings}}
\author{Kazunori Hayashi\thanks{Kyoto University:~~\texttt{kazunori@i.kyoto-u.ac.jp}}
\and
Yasuaki Hiraoka\thanks{Hiroshima University/JST:~~\texttt{hiraok@hiroshima-u.ac.jp}}
}

\begin{document}
\maketitle
\begin{abstract}
This paper studies maximum likelihood(ML) decoding in error-correcting codes as rational maps
and proposes an approximate ML decoding rule by using a Taylor expansion.
The point for the Taylor expansion, which will be denoted by $p$ in the paper, 
is properly chosen by considering some dynamical system properties.
We have two results about this approximate ML decoding.
The first result proves that the order of the first nonlinear terms in the Taylor expansion is determined 
by the minimum distance of its dual code. 
As the second result, we give numerical results on bit error probabilities for the approximate ML decoding. 
These numerical results show better performance than that of BCH codes, and 
indicate that this proposed method approximates the original ML decoding very well.
\end{abstract}
{\small {\bf Key words.} Maximum likelihood decoding, rational map, dynamical system\vspace{0.1cm}\\
{\bf AMS subject classification.} 37N99, 94B35}
\section{Introduction}\label{sec:intro}
This paper proposes a new perspective to maximum likelihood(ML) decoding 
in error-correcting codes as rational maps and 
shows some relationships between coding theory and dynamical systems.
In Section \ref{sec:cs}, \ref{sec:lc}, and \ref{sec:ml}, 
we explain notations and minimum prerequisites of coding theory (e.g., see \cite{Gallagerinfo}).
The main results are presented in Section \ref{sec:result}.

\subsection{Communication Systems}\label{sec:cs}

\begin{figure}[htbp]
\begin{minipage}{0.5\hsize}
\begin{center}
\includegraphics[width=8.0cm,height=2.5cm]{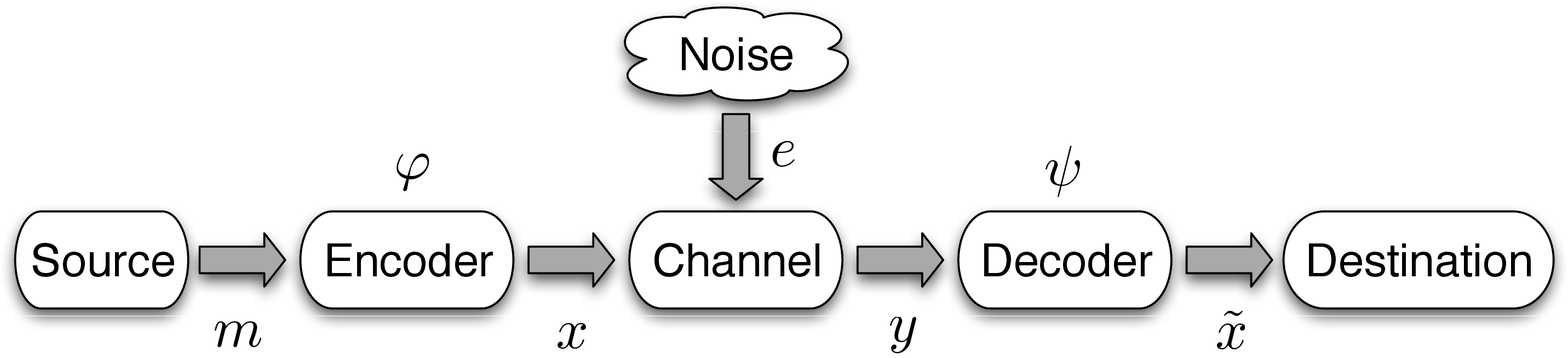}
\end{center}
\caption{Communication system}
\label{fig:channel}
\end{minipage}
\begin{minipage}{0.5\hsize}
\begin{center}
\includegraphics[width=4cm,height=3.0cm]{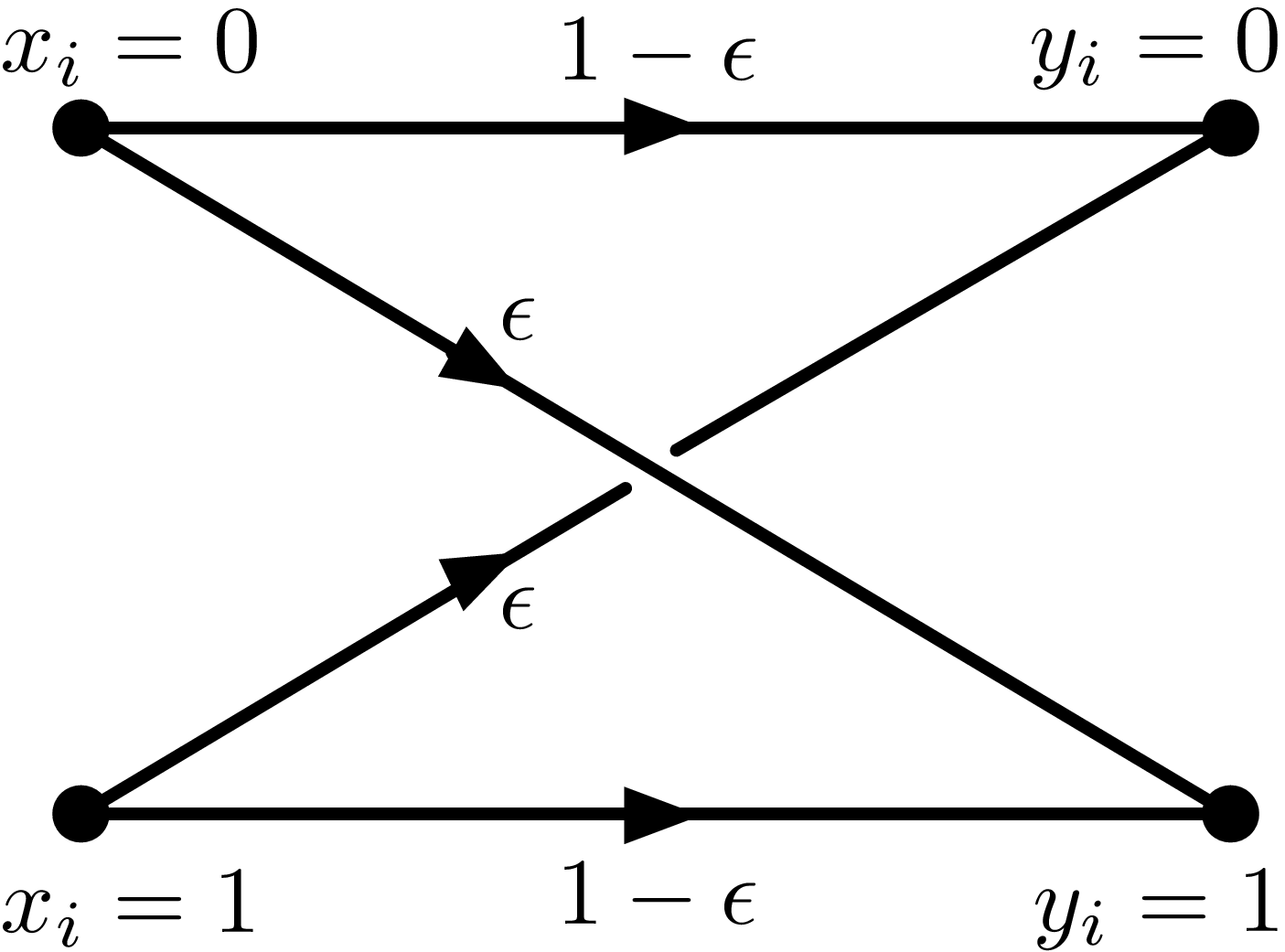}
\end{center}
\caption{Binary symmetric channel}
\label{fig:bsc}
\end{minipage}
\end{figure}

A mathematical model of communication systems in information theory was developed by
Shannon \cite{shannon}.  A general block diagram for visualizing the behavior of such systems is given by
Figure \ref{fig:channel}. 
The source transmits a $k$-bit message $m=(m_1\cdots m_k)$ to the destination via the channel, 
which is usually affected by noise $e$. 
In order to recover the transmitted message at the destination under the influence of 
noise, we transform the message into a codeword $x=(x_1\cdots x_n),~n\geq k,$ by some injective mapping $\varphi$
at the encoder and input it to the channel. 
Then the decoder transforms an $n$-bit received sequence of letters $y=(y_1\cdots y_n)$ by
some decoding mapping $\psi$ in order to obtain the transmitted codeword at the destination. 
Here we consider all arithmetic calculations in some finite field and in this paper we fix it as $\mathbb{F}_2=\{0,1\}$ 
except for Section \ref{sec:ag}.
As a model of channels, we deal with a binary symmetric channel (BSC) in this paper which is characterized by the 
transition probability $\epsilon$ ($0<\epsilon<1/2$). Namely, with probability $1-\epsilon$, the output letter is a faithful replica of the 
input, and with probability $\epsilon$, it is the opposite of the input letter for each bit (see Figure \ref{fig:bsc}). 
In particular, this is an example of memoryless channels.

Then, one of the main purposes of coding theory is to develop a good encoder-decoder pair $(\varphi, \psi)$ 
which is robust to noise perturbations. 
Hence, the problem is how we efficiently use the redundancy $n\geq k$ in this setting.
\subsection{Linear Codes}\label{sec:lc}
A code with a linear encoding map $\varphi$ is called a linear code. A codeword in a linear code can be 
characterized by its generator matrix
\[
G=\left(\begin{array}{ccc}
g_{11} & \cdots & g_{1n}\\
\vdots & \ddots & \vdots\\
g_{k1} & \cdots & g_{kn}
\end{array}
\right)=
\left(
g_1 \cdots g_n
\right),~~~
g_i=\left(\begin{array}{c}
g_{1i}\\
\vdots \\
g_{ki}
\end{array}
\right),~~i=1,\cdots,n,
\]
where each element $g_{ij}\in\mathbb{F}_2$.
Therefore the set of codewords is given by
\[
\mathcal{C}=\{
(x_1\cdots x_n)=(m_1\cdots m_k)G~|~m_i\in \mathbb{F}_2,\},
~~~\sharp \mathcal{C}=2^k.
\]
Here without loss of generality, we assume $g_i\neq 0$ for all $i=1,\cdots, n$ and ${\rm rank}~G=k$. 
We call $k$ and $n$ the dimension and the length of the code, respectively.

Because of the linearity, it is also possible to describe $\mathcal{C}$ as a kernel of a matrix $H$ 
whose $m=n-k$ row vectors are linearly independent and orthogonal to those of $G$, i.e.,
\[
\mathcal{C}=\{(x_1\cdots x_n)~|~(x_1\cdots x_n)H^t=0\},~~~
H=\left(\begin{array}{ccc}
h_{11} & \cdots & h_{1n}\\
\vdots & \ddots & \vdots\\
h_{m1} & \cdots & h_{mn}
\end{array}
\right)=
(h_1\cdots h_n),
\]
where $H^t$ means the transpose matrix of $H$. This matrix $H$ is called a parity check matrix of $\mathcal{C}$.
The dual code $\mathcal{C}^*$ of $\mathcal{C}$ is defined in such a way that a parity check matrix of 
$\mathcal{C}^*$ is given by a generator matrix $G$ of $\mathcal{C}$.

The Hamming distance $d(x,y)$ between two $n$-bit sequences $x,y\in\mathbb{F}^n_2$ is given
by the number of positions at which the two sequences differ. The weight of an element $x\in\mathbb{F}^n_2$ is
the Hamming distance to $0$, i.e., $d(x):=d(x,0)$.
Then the minimum distance $d(\mathcal{C})$ of a code $\mathcal{C}$ is defined by two different ways as
\[
d(\mathcal{C})=\min \{
d(x,y)~|~x,y\in\mathcal{C}~{\rm and}~x\neq y
\}
=\min \{
d(x)~|~0\neq x\in\mathcal{C}
\}.
\]
Here the second equality results from the linearity. It is easy to observe that the minimum distance is $d(\mathcal{C})=d$
if and only if there exists a set of $d$ linearly dependent column vectors of $H$ but no set of $d-1$ linearly dependent 
column vectors. 

For a code $\mathcal{C}$ with the minimum distance $d=d(\mathcal{C})$, let us set $t:=\lfloor (d-1)/2\rfloor$, where
$\lfloor a \rfloor$ is the integer part of $a$. 
Then, it follows from the following observation that $\mathcal{C}$ can correct $t$ errors: 
if $y\in\mathbb{F}^n_2$ and $d(x,y)\leq t$ for some $x\in\mathcal{C}$ then $x$ is the 
only codeword with $d(x,y)\leq t$. 
In this sense, the minimum distance is one of the important parameters to measure performance of a code and
is desirable to design it as large as possible for the robustness to noise.
%
%
%
%
%
\subsection{Maximum Likelihood Decoding}\label{sec:ml}
Let us recall that, given a transmitted codeword $x$, the conditional probability $P(y|x)$ of a received sequence 
$y\in\mathbb{F}^n_2$ at the decoder is given by
\[
P(y|x)=P(y_1|x_1)\cdots P(y_n|x_n)
\]
for a memoryless channel. 
Maximum likelihood(ML) decoding $\psi:\mathbb{F}_2^n\ni y\mapsto \tilde{x}\in\mathbb{F}^n_2$ is given by taking 
the marginalization of $P(y|x)$ for each bit. 
Precisely speaking, for a received sequence $y$, the $i$-th bit element $\tilde{x}_i$ of the decoded word $\tilde{x}$ is determined by the following rule:
\begin{equation}
\tilde{x}_i:=\left\{
\begin{array}{ll}
1, & \sum_{\substack{x\in\mathcal{C}\\x_i=0}}P(y|x)\leq\sum_{\substack{x\in\mathcal{C}\\x_i=1}}P(y|x)\\
0, & {\rm otherwise}
\end{array},~~i=1,\cdots,n.
\right. \label{eq:bitml}
\end{equation}

In general, for a given decoder $\psi$, the bit error probability ${P^e_{\rm bit}=\max\{P^e_1,\cdots,P^e_n\}}$,
where
\[
P^e_i=\sum_{x\in\mathcal{C},y\in\mathbb{F}^n_2}P(x,y)(1-\delta_{x_i,\tilde{x}_i}),~~~
\tilde{x}=(\tilde{x}_1,\cdots,\tilde{x}_n)=\psi(y),~~~
\delta_{a,b}=\left\{
\begin{array}{ll}
1, & a=b\\
0, & a\neq b
\end{array}
\right.,
\]
is one of the important measures of decoding performance. 
Obviously, it is desirable to design an encoding-decoding pair whose bit error probability is as small as possible.
It is known that ML decoding attains the minimum bit error probability for any encodings 
under the uniform distribution on $P(x)$. 
In this sense, ML decoding is the best for all decoding rules. 
However its computational cost requires at least $2^k$ operations, and it is too much to use for practical applications.

From the above property of ML decoding, 
one of the key motivation of this work comes from the following simple question. 
Is it possible to accurately approximate the ML decoding rules with low computational complexity? 
The main results in this paper give answers to this question.

\subsection{Main Results}\label{sec:result}

Let us first define for each codeword $x\in \mathcal{C}$ its codeword polynomial $F_x(u_1,\cdots,u_n)$ as
\[
F_x(u_1,\cdots,u_n):=\prod^n_{i=1}\rho_i(u_i), ~~~~~
\rho_i(u_i)=\left\{\begin{array}{ll}
u_i,& x_i=1,\\
1-u_i,&x_i=0.
\end{array}
\right.
\]
Then we define a rational map $f:I^n\rightarrow I^n$, $I=[0,1],$ by using codeword polynomials as
\begin{eqnarray}
&&f:(u_1,\cdots,u_n)\mapsto (u'_1,\cdots,u'_n),\nonumber\\
&&u'_i=f_i(u):=\frac{\sum_{x\in \mathcal{C},x_i=1} F_x(u)}{H(u)},~~~i=1,\cdots,n,\nonumber\\
&&H(u):=\sum_{x\in\mathcal{C}}F_x(u),\label{eq:map}
\end{eqnarray}
where $u=(u_1,\cdots,u_n)$.
This rational map plays the most important role in the paper. It is sometimes denoted 
by $f_G$, when we need to emphasize the generator matrix $G$ of the code $\mathcal{C}$. 

For a sequence $y\in\mathbb{F}^n_2$, let us take a point $u^0\in I^n$ as
\begin{equation}
u^0_i=\left\{
\begin{array}{ll}
\epsilon, & y_i=0\\
1-\epsilon, & y_i=1
\end{array}
,~~~i=1,\cdots,n, \label{eq:ip}
\right.
\end{equation}
where $\epsilon$ is the transition probability of the channel. 
Then it is straightforward to check that $F_x(u^0)=P(y|x)$. Namely, the conditional probability of $y$ under a codeword 
$x\in \mathcal{C}$ is given by the value of the corresponding codeword polynomial $F_x(u)$ at $u=u^0$.
Therefore, from the construction of the rational map, ML decoding (\ref{eq:bitml}) is equivalently given 
by the following rule
\begin{eqnarray}
&&\psi:\mathbb{F}_2^n\ni y\mapsto \tilde{x}\in\mathbb{F}_2^n,\nonumber\\
&&\tilde{x}_i=\psi_i(y):=\left\{
\begin{array}{ll}
1, & f_i(u^0)\geq 1/2\\
0, & f_i(u^0)< 1/2
\end{array}
,~~~i=1,\cdots,n.
\right.\label{eq:mlrule}
\end{eqnarray}
In this sense, the study of ML decoding can be treated by analyzing the image of the initial 
point (\ref{eq:ip}) by the rational map (\ref{eq:map}).
Some of the properties of this map in the sense of dynamical systems will be studied in detail in Section \ref{sec:ds}. 
We will also discuss in Section \ref{sec:discussions} that performance of a code can be explained by these properties.

For the statement of the main results, we only here mention that this rational map has a fixed point 
$p:=(1/2,\cdots,1/2)$ for any generator matrix (Proposition \ref{pp:fp}). 
Let us denote the Taylor expansion at $p$ by
\begin{equation}
f(u)=p+Jv+f^{(2)}(v)+\cdots+f^{(l)}(v)+O(v^{l+1}), \label{eq:taylor}
\end{equation}
where $v=(v_1,\cdots,v_n)$ is a vector notation of $v_i=u_i-1/2, i=1,\cdots,n$, 
$J$ is the Jacobi matrix at $p$, $f^{(i)}(v)$ corresponds to the $i$-th order term, and 
$O(v^{l+1})$ means the usual order notation. 
The reason why we choose $p$ as the approximating point is related to the local dynamical property at $p$ 
and will be explained in Section \ref{sec:p}.

By truncating higher oder terms $O(v^{l+1})$ in (\ref{eq:taylor}) and denoting it as 
\[
\tilde{f}(u)=p+Jv+f^{(2)}(v)+\cdots+f^{(l)}(v), 
\]
we can define the $l$-th approximation of ML decoding by replacing the map $f(u)$ in 
(\ref{eq:mlrule}) with $\tilde{f}(u)$, and 
denote this approximate ML decoding by ${\tilde{\psi}:\mathbb{F}_2^n\rightarrow\mathbb{F}_2^n}$, i.e.,
\begin{eqnarray}
&&\tilde\psi:\mathbb{F}_2^n\ni y\mapsto \tilde{x}\in\mathbb{F}_2^n,\nonumber\\
&&\tilde{x}_i=\tilde\psi_i(y):=\left\{
\begin{array}{ll}
1, & \tilde{f}_i(u^0)\geq 1/2\\
0, & \tilde{f}_i(u^0)< 1/2
\end{array}
,~~~i=1,\cdots,n.
\right.\label{eq:apmlrule}
\end{eqnarray}
Let us remark that the notations $\tilde{f}$ and $\tilde{\psi}$ do not explicitly express the dependence on $l$
for removing unnecessary confusions of subscripts.

\subsubsection{Duality Theorem}\label{sec:duality}
We note that there are two different viewpoints on this approximate ML decoding.
One way is that, in the sense of its precision, it is preferable to have an expansion with large $l$. 
On the other hand, from the viewpoint of low computational complexity, it is desirable to include many zero elements in 
higher order terms. 
The next theorem states a sufficient condition to satisfy these two requirements. 
\begin{thm}\label{thm:thm}
Let $l\geq 2$. If any distinct $l$ column vectors of a generator matrix $G$ are linearly independent, 
then the Taylor expansion {\rm (\ref{eq:taylor})} at $p$ of the rational map {\rm (\ref{eq:map})} 
takes the following form
\[
f(u)=u+f^{(l)}(v)+O(v^{l+1}),
\]
where the $i$-th coordinate $f^{(l)}_i(v)$ of $f^{(l)}(v)$ is given by
\begin{eqnarray}
&&f^{(l)}_i(v)=\sum_{(i_1,\cdots,i_l)\in\Theta^{(l)}_i}(-2)^{l-1}
v_{i_1}\cdots v_{i_l}=
-\frac{1}{2}\!\!\!\!\sum_{(i_1,\cdots,i_l)\in\Theta^{(l)}_i}(1-2u_{i_1})\cdots(1-2u_{i_l}),\label{eq:hosum}\\
&&\Theta^{(l)}_i=\left\{(i_1,\cdots,i_l)~|~1\leq i_1<\cdots <i_l\leq n,~ 
i_k\neq i~(k=1,\cdots,l),~
g_{i}+g_{i_1}+\cdots+g_{i_l}=0
\right\}.\nonumber
\end{eqnarray}
\vspace{-0.2cm}
\end{thm}

First of all, it follows that the larger the minimum distance of the dual code 
$\mathcal{C}^*$ is, the more precise approximation of ML decoding with low computational complexity 
we have for the code $\mathcal{C}$ with the generator matrix $G$. 
Especially, we can take $l=d(\mathcal{C^*})-1$.

Secondly, let us consider the meaning of the approximate map $\tilde{f}(u)$ and 
its approximate ML decoding $\tilde{\psi}$. 
We note that each value $u^0_i(i=1,\cdots,n)$ in (\ref{eq:ip})
for a received word $y\in\mathbb{F}^n_2$ expresses the likelihood $P(y_i|x_i=1)$. 
Let us suppose $\Theta^{(l)}_i\neq\emptyset$. 
Then, from the definition of $u^0$, each term in the sum of (\ref{eq:hosum}) satisfies
\[
-\frac{1}{2}(1-2u_{i_1})\cdots(1-2u_{i_l})
\left\{
\begin{array}{ll}
<0, & {\rm if}~y_{i_1}+\cdots+y_{i_l}=0\\
>0, & {\rm if}~y_{i_1}+\cdots+y_{i_l}=1
\end{array},
\right.
~~~~(i_1,\cdots,i_l)\in\Theta^{(l)}_i.
\]
When $y_{i_1}+\cdots+y_{i_l}=0(=1, {\rm resp.})$, this term decreases(increases, resp.) the value of initial likelihood $u^0_i$.
In view of the decoding rule (\ref{eq:mlrule}), this induces $\tilde{x}_i$ to be decoded 
into $\tilde{x}_i=0(=1,{\rm resp.})$, 
and this actually corresponds to the structure of the code $g_{i}+g_{i_1}+\cdots+g_{i_l}=0$
appearing in $\Theta^{(l)}_i$. 
In this sense, the approximate map ${\tilde f}(u)$ can be regarded as renewing the likelihood 
(under suitable normalizations) based on the code structure,
and the approximate ML decoding $\tilde{\psi}$ judges these renewed data.
From this argument, it is easy to see that a received word $y\in\mathcal{C}$ is decoded into $y=\tilde{\psi}(y)\in\mathcal{C}$,
i.e., the codeword is decoded into itself and, of course, this property should be equipped with any decoders. 

We also remark that Theorem \ref{thm:thm} can be regarded as a duality theorem in the following sense.
Let $\mathcal{C}$ be a code whose generator(resp. parity check) matrix is $G$ (resp. $H$).
As we explained in Section \ref{sec:lc}, the linear independence of the column vectors of $H$ controls
the minimum distance $d(\mathcal{C})$ and this is an encoding property.
On the other hand, Theorem \ref{thm:thm} shows that the linear independence of the column vectors of $G$, 
which determines the dual minimum distance $d(\mathcal{C}^*)$, controls 
a decoding property of ML decoding in the sense of accuracy and computational complexity.
Hence, we have the correspondence between $H/G$ duality and encoding/decoding duality. 
In Corollary \ref{corollary:ag}, we will consider this duality viewpoint 
in a setting of geometric Reed-Solomon/Goppa codes.


\subsubsection{Decoding Performance}\label{sec:dprc}
\begin{figure}[htbp]
\begin{minipage}{0.5\hsize}
\begin{center}
\includegraphics[width=8.0cm]{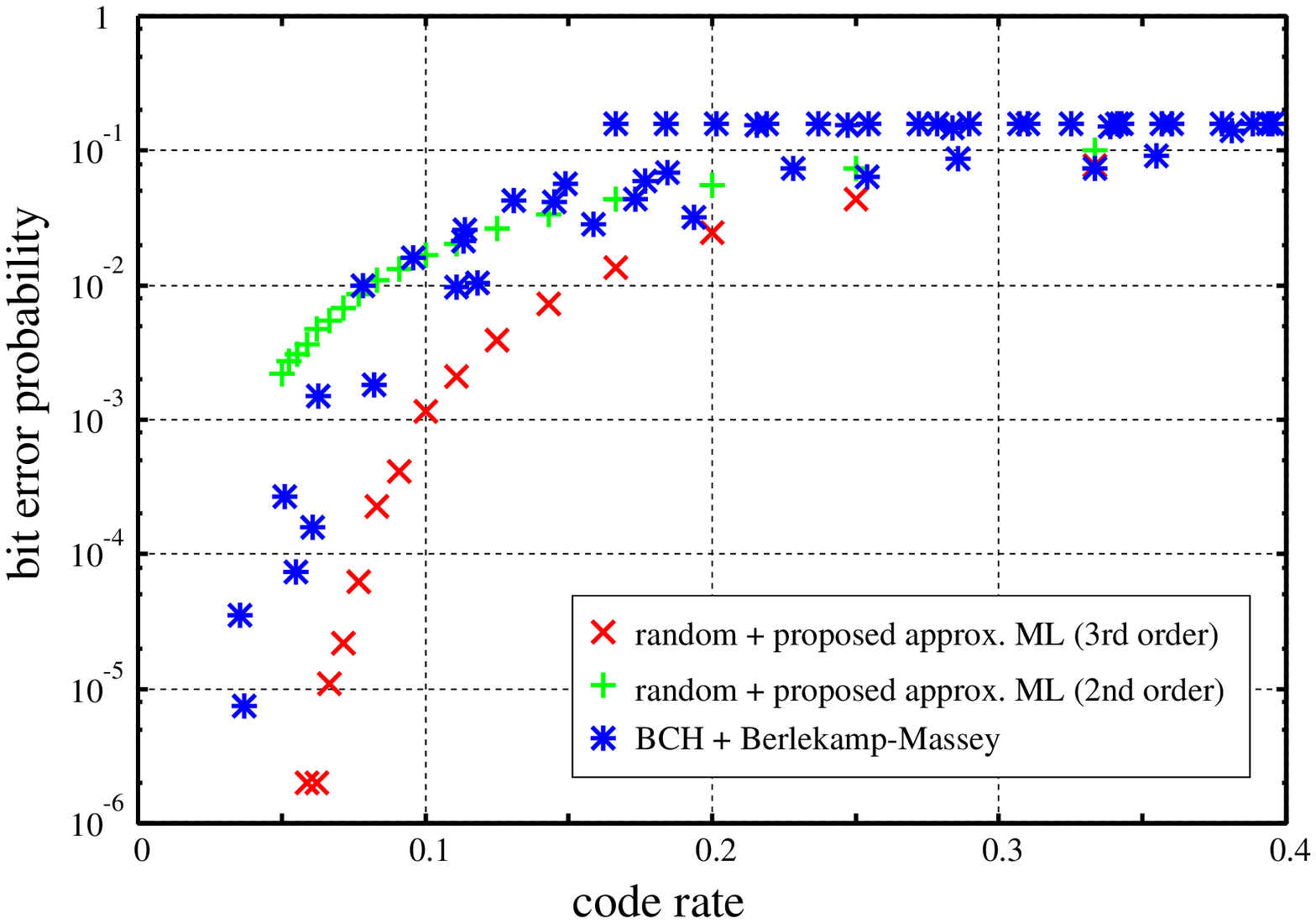}\vspace{-0.5cm}
\end{center}
\caption{BSC with $\epsilon=0.16$. 
Horizontal axis: code rate ${r=k/n}$.
Vertical axis: bit error probability. 
$\times$: random codes with the 3rd order approximate ML decoding.  
$+$: random codes with the 2nd order approximate ML decoding.
$\ast$: BCH codes with Berlekamp-Massey decoding.
}
\label{fig:berrate}
\end{minipage}\hspace{0.5cm}
\begin{minipage}{0.5\hsize}
\begin{center}
\includegraphics[width=8.0cm]{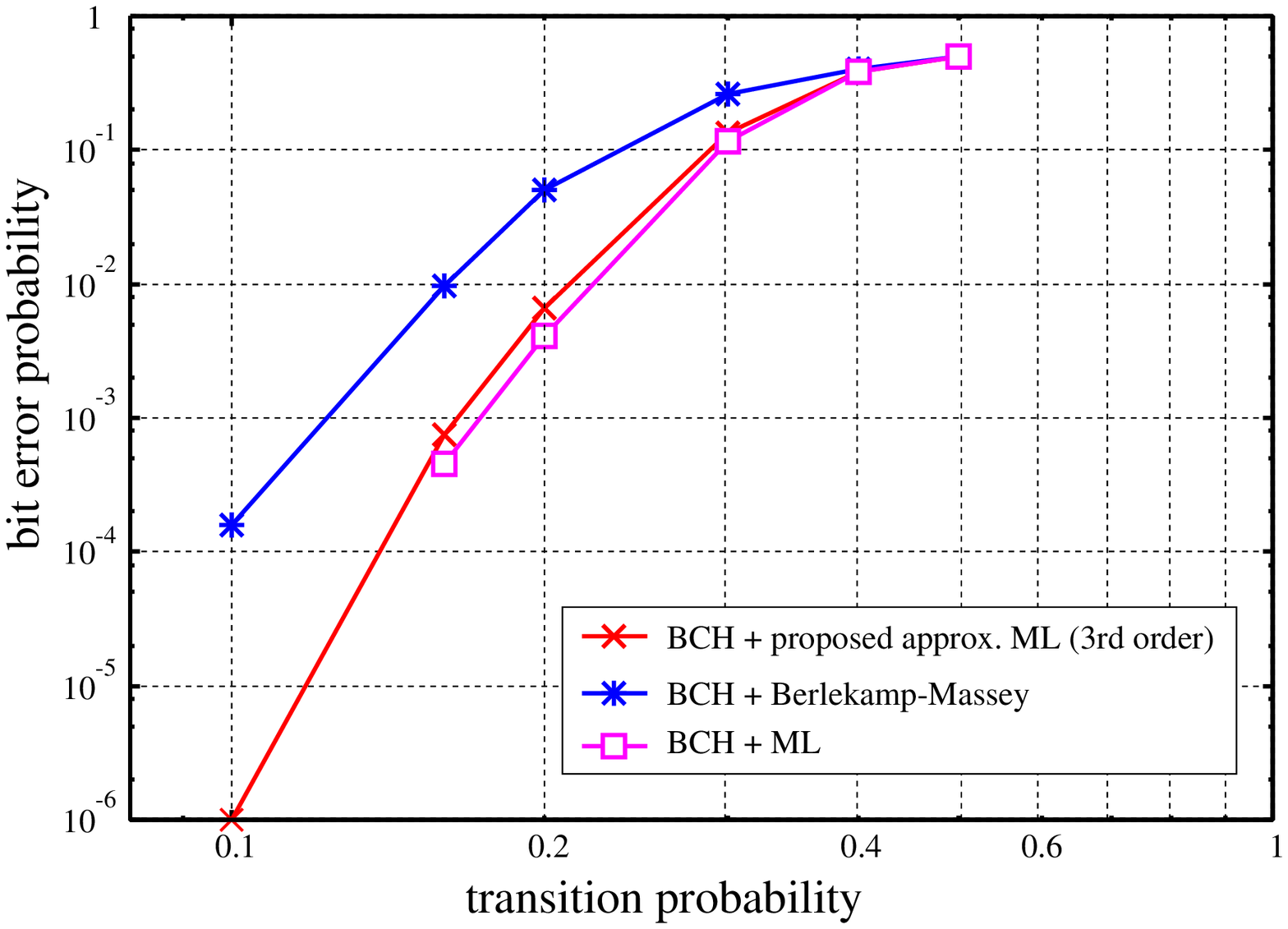}\vspace{-0.5cm}
\end{center}
\caption{
Comparison of decoding performances for a BCH code with $k=7$ and $n=63$. 
Horizontal axis: transition probability $\epsilon$. 
Vertical axis: bit error probability. 
$\times$: 3rd order approximate ML decoding.
$\ast$: Berlekamp-Massaey decoding. 
$\Box$: ML decoding.
}
\label{fig:berprob}
\end{minipage}
\end{figure}

We show the second result of this paper about the decoding performances of the approximate ML (\ref{eq:apmlrule}). 
For this purpose, let us first examine numerical simulations of the bit error probability on the BSC with 
the transition probability $\epsilon=0.16$.
We also show numerical results on BCH codes with Berlekamp-Massey decoding for comparison.
The results are summarized in Figure \ref{fig:berrate}.

Here, the horizontal axis is the code rate $r=k/n$, and the vertical axis is the bit error probability.
The plots $+$ ($\times$ resp.) correspond to the 2nd (3rd resp.) order approximate ML (\ref{eq:apmlrule}), 
and $\ast$ are the results on several BCH codes ($n=7,~15,~31,~63,~127,~255,~511$) with Berlekamp-Messey decodings.
For the proposed method,
we randomly construct a systematic generator matrix in such a way that each column except
for the systematic part has the same weight (i.e. number of non-zero elements) $w$.
To be more specific, the submatrix composed by the first $k$ columns of the generator matrix is set to be 
an identity matrix in order to make the code systematic, while the rest of the generator matrix is made up
of $k\times k$ random matrices generated by random permutations of columns of a circulant matrix, 
whose first column is given by
\begin{align}
(\underbrace{1,\cdots 1}_{w},0,\cdots 0)^{t} \notag.
\end{align}
The reason for using random codings is that we want to investigate average behaviors of the decoding performance,
and, for this purpose, we do not put unnecessary additional structure at encodings.
The number of matrices added after the systematic part depends on the code rate, and the plot for each code rate
corresponds to the best result obtained out of about 100 realizations of the generator matrix.
Also, we have employed $w=2$ and 3 for the generator matrices of the 3rd and the 2nd order approximate ML, respectively.
Moreover, the length of the codewords $n$ are assumed to be up to 512.
From Figure \ref{fig:berrate}, we can see that the proposed method with the 3rd order approximate ML ($\times$) 
achieves better performance than that of BCH codes with Berlekamp-Massey($\ast$).
It should be also noticed that the decoding performance is improved a lot from the 2nd order to the 3rd order approximation. 
This improvement is reasonable because of the meaning of the Taylor expansion.
 
Next, let us directly compare the decoding performances among ML, approximate ML (3rd order) and Berlekamp-Massey by
applying them on the same BCH code ($k=7$, $n=63$).
The result on the bit error probability with respect to transition probability is shown in Figure \ref{fig:berprob}.
This figure clearly shows that the performance of the 3rd order approximate ML decoding is far better than that of 
Berlekamp-Massey decoding (e.g., improvement of double-digit at $\epsilon=0.1$). 
Furthermore, it should be noted that the 3rd order approximate ML decoding achieves a very close bit error performance
to that of ML decoding. 
Although we have not mathematically confirmed the computational complexity of the proposed approach,
the computational time of the approximate ML (3rd order) is much faster than ML decoding.
This fact about low computational complexity of the approximate ML is explained as follows: 
Non-zero higher order terms in (\ref{eq:hosum}) appear as a result of linear dependent relations of 
column vectors of $G$, however, linear dependences require high codimentional scenario. 
Hence, most of the higher order terms become zeros. 
As a result, the computational complexity for the approximate ML, which is determined by the number 
of nonzero terms in the expansion, becomes small.

In conclusion, these numerical simulations suggest that the 3rd order approximate ML decoding approximates 
ML decoding very well with low computational complexity.
We notice that the encodings examined here are random codings.
Hence, we can expect to obtain better bit error performance by introducing certain structure on encodings suitable 
to this proposed decoding rule, or much more exhaustive search of random codes.
One of the possibility will be the combination with Theorem \ref{thm:thm}. 
On the other hand, it is also possible to consider suitable encoding rules from the viewpoint of 
dynamical systems via rational maps (\ref{eq:map}). 
This issue is discussed in detail in Section \ref{sec:discussions}.
In any case, finding suitable encoding structure for the proposed decodings is one of the important future problem.

\vspace{0.5cm}
The paper is organized as follows. In Section \ref{sec:ds}, we study properties of the rational map (\ref{eq:map}) 
in view of a discrete dynamical system and show some relationships to coding theory. 
We also show that this discrete dynamical system is related to 
a continuous gradient dynamical system $du/dt={\rm grad}(\log H(u))$.
Section \ref{sec:taylor} deals with relationship between a generator matrix of a code and its Taylor expansion 
(\ref{eq:taylor}). The proof of Theorem \ref{thm:thm}, which is a direct consequence of Proposition \ref{pp:jacobi} and 
Proposition \ref{pp:fl0}, is shown in this section. 
In Section \ref{sec:ag}, we apply Theorem \ref{thm:thm} to geometric Reed-Solomon/Goppa codes in Corollary \ref{corollary:ag}
with a simple example by using the Hermitian curve.
Finally, we discuss the future problems on this subject as an intersection of dynamical systems and coding theory.
\section{Rational Maps}\label{sec:ds}
In this section, we discuss ML decoding from dynamical systems viewpoints.
Let $G$ be a $k\times n$ generator matrix for a code $\mathcal{C}$. We begin with showing the following easy
consequence of linear codes, which will be used frequently throughout the paper.
\begin{lemma}\label{lemma:gi}
For $i\in\{1,\cdots,n\}$, let us denote subcodes of $\mathcal{C}$ with $x_i=0$ and $x_i=1$ respectively by
\[
\mathcal{C}(x_i=0)=\left\{x\in \mathcal{C}~|~x_i=0\right\},~~~~~
\mathcal{C}(x_i=1)=\left\{x\in \mathcal{C}~|~x_i=1\right\}.
\]
Then $\sharp \mathcal{C}(x_i=0)=\sharp \mathcal{C}(x_i=1)=2^{k-1}$.
\end{lemma}
\begin{prf}{\rm 
From the assumption on the generator matrix, we have $g_i\neq 0$ and let $l>0$ be the number of 1 in $g_i$.
Then we have
\[
\sharp \mathcal{C}(x_i=0)=
\sum^l_{\substack{
j=0 \\
j:{\rm even}
}}{_l}C_j
\times 2^{k-l},~~~~
\sharp \mathcal{C}(x_i=1)=
\sum^l_{\substack{
j=0 \\
j:{\rm odd}
}}{_l}C_j
\times 2^{k-l},
\]
where the symbol ${_l}C_j$ means the number of combinations for taking $j$ elements from $l$ elements.
However these summations of combinations are obviously equal because
\[
0=(1-1)^l=\sum^l_{j=0}{_l}C_j(-1)^j=
\sum^l_{\substack{
j=0 \\
j:{\rm even}
}}{_l}C_j
-
\sum^l_{\substack{
j=0 \\
j:{\rm odd}
}}{_l}C_j.
\]
Therefore $\sharp \mathcal{C}(x_i=0)=\sharp \mathcal{C}(x_i=1)=2^{k-1}$.
}\qed
\end{prf}
Next, let us characterize the codewords in $\mathcal{C}$ and the non-codewords in $\mathbb{F}^n_2\setminus \mathcal{C}$ 
by means of the rational map (\ref{eq:map}).
It should be remarked that a codeword polynomial $F_x(u)$ 
for $u=(u_1,\cdots,u_n)$ with $u_i\in\partial I=\{0,1\},i=1,\cdots,n$, 
takes its value
\[
F_x(u)=\left\{\begin{array}{ll}
1, & u=x,\\
0, & u\neq x.
\end{array}
\right.
\]
Here we identify a point $u=(u_1,\cdots,u_n)\in I^n,$ ${u_i\in\partial I, i=1,\cdots,n}$,
with a point $u\in\mathbb{F}^n_2$ by a natural inclusion $\mathbb{F}_2\hookrightarrow I$, 
and this convention will be used frequently in the paper.
We also define a point $u\in I^n$ as a pole of the rational map (\ref{eq:map}) if $H(u)=0$.
The boundary and the interior of a set $A$ are denoted by $\partial A$  and ${\rm Int}(A)$, respectively.

\begin{pp}\label{pp:fp}
The followings hold for the rational map {\rm (\ref{eq:map})}:
\begin{enumerate}
\item $p:=(1/2,\cdots,1/2)\in I^n$ is a fixed point.
\item Let $u\in\mathbb{F}^n_2$. 
Then $u$ is a fixed point if and only if $u\in\mathcal{C}$.
\item The set of poles is given by
\[
\hspace{-0.5cm}S:=\left\{
u\in \partial (I^n)~
\begin{array}{|l}
there~is~no~code~word~x\in\mathcal{C}~with~(x_{i_1},\cdots,x_{i_l})=(u_{i_1},\cdots,u_{i_l}),\\
where~u_{i_1},\cdots,u_{i_l}\in \mathbb{F}_2,~and~0<u_i<1,~i\neq i_1,\cdots,i_l 
\end{array}
\right\}.
\]
Especially, $\mathbb{F}^n_2\setminus\mathcal{C}\subset S$.
\end{enumerate}
\end{pp}
\begin{prf}{\rm 
For the statement 1, let us note that $F_x(p)=(1/2)^n$ for each codeword $x\in\mathcal{C}$.
Then Lemma \ref{lemma:gi} leads to $f_i(p)=1/2, i=1,\cdots,n.$

For the statement 2, let us suppose $u\in\mathcal{C}$. Then, from the remark before this proposition,
$f_i(u)=1({\rm or}~0~{\rm resp.})$ if $u_i=1({\rm or}~0~{\rm resp.})$ for each $i$. Hence $f(u)=u$. 
On the other hand, if $u\in\mathbb{F}^n_2\setminus\mathcal{C}$, then $H(u)=0$. It means that $u$ is a pole and 
can not be a fixed point.

For the statement 3, let us first note that $u\in S$ if and only if $F_x(u)=0$ for any codeword $x\in\mathcal{C}$, 
because $F_x(u)\geq 0$ for $u\in I^n$. Therefore $S\subset \partial (I^n)$, since $F_x(u)>0$ for $u\in{\rm Int}(I^n)$.

Let $u\in\partial (I^n)$ such that $u_{i_1},\cdots,u_{i_l}\in\mathbb{F}_2$ and $0<u_i<1$ for $i\neq i_1,\cdots,i_l$.
Then, if there exists a codeword $x\in\mathcal{C}$ such that $(x_{i_1},\cdots,x_{i_l})=(u_{i_1},\cdots,u_{i_l})$, 
then the value of its corresponding codeword polynomial is $F_x(u)>0$. So, $u\notin S$.
On the other hand, if there is no such codeword, then $H(u)=0$, it means $u\in S$.
}\qed
\end{prf}

From this proposition, the rational map (\ref{eq:map}) has information of not only all codewords $\mathcal{C}$ 
as fixed points but also non-codewords $\mathbb{F}^n_2\setminus\mathcal{C}$ as poles. 
We call these fixed points codeword fixed points. 
The following proposition shows that all of the codeword fixed points are stable. 
\begin{pp}\label{pp:stable}
Let a parity check matrix do not have zero column vectors, i.e., there exists no codeword with weight 1.
Let $u$ be a codeword fixed point. Then the Jacobi matrix of the rational map {\rm (\ref{eq:map})} at $u$ is the zero matrix. 
Hence, $u$ is a stable fixed point.
\end{pp}
\begin{prf}{\rm 
Let us denote the $i$-th element of the rational map (\ref{eq:map}) by
\[
f_i(u)=\frac{I_i(u)}{H(u)}
\]
and denote the derivatives of $I_i$ and $H$ with respect to $u_j$ by $I^j_i$ and $H^j$, respectively, for the simplicity of notations.
In what follows, we will also use these notations for higher order derivatives in a similar way 
(like $H^{ij}=\frac{\partial^2 H}{\partial u_i \partial u_j}$).
Then the derivative $\partial f_i/\partial u_j(u)$ is given by
\[
\frac{\partial f_i}{\partial u_j}(u)=\left\{
I^j_i(u)H(u)-I_i(u)H^j(u)
\right\}/H^2(u).
\]
Since $u\in\mathcal{C}$, we have $H(u)=1$. Let us consider the two cases $u_i=0$ and $u_i=1$, separately.

For $u_i=0$, we have $I_i(u)=0$. Similarly we have $I^j_i(u)=0$ if $i\neq j$.
Hence, $\partial f_i/\partial u_j(u)=0$ in this case. 
On the other hand, if $i=j$, from the assumption on the parity check matrix, 
we have no codeword $\tilde{u}\in\mathcal{C}$ such that the only difference from $u$ 
occurs at the $i$-th bit element, i.e., $u_j=\tilde{u}_j$ for $j\neq i$ but $u_i\neq \tilde{u}_i$. 
Therefore, $I^j_i(u)=0$, and it again leads to $\partial f_i/\partial u_j(u)=0$.
The case $u_i=1$ can be proven by the similar way.
}\qed
\end{prf}

Let us next discuss properties of the fixed point $p=(1/2,\cdots,1/2)$. 
\begin{pp}\label{pp:jacobi}
Let $G=(g_1 \cdots g_n)$ be a generator matrix. Then, the Jacobi matrix $J$ of the rational map {\rm (\ref{eq:map})} 
at $p$ is determined by
\[
\begin{array}{ll}
J_{ij}=1,  & {\rm if}~g_i=g_j,\\
J_{ij}=0,  & {\rm if}~g_i\neq g_j\\
\end{array}
\]
for all $i,j=1,\cdots,n$. 
\end{pp}
For the proof of this proposition, we need the following lemma.
\begin{lemma}\label{lemma:gij}
For $i,j\in\{1,\cdots,n\}$ with $i\neq j$, let us consider the following subcodes of $\mathcal{C}$
\[
\mathcal{C}(\substack{x_i=1\\x_j=0})=\left\{
x\in\mathcal{C}~|~x_i=1,x_j=0
\right\},~~~~
\mathcal{C}(\substack{x_i=1\\x_j=1})=\left\{
x\in\mathcal{C}~|~x_i=1,x_j=1
\right\}.
\]
Then the followings hold
\begin{enumerate}
\item if $g_i=g_j$, then $\sharp \mathcal{C}(\substack{x_i=1\\x_j=0})=0,~\sharp \mathcal{C}(\substack{x_i=1\\x_j=1})=2^{k-1}.$
\item if $g_i\neq g_j$, then $\sharp \mathcal{C}(\substack{x_i=1\\x_j=0})=\sharp \mathcal{C}(\substack{x_i=1\\x_j=1})=2^{k-2}.$
\end{enumerate}
\end{lemma}
\begin{prf}{\rm 
The statement is trivial when $g_i=g_j$, so we suppose $g_i\neq g_j$. 
In the following we adopt the $mod~2$ arithmetic for elements in $\mathbb{F}_2$.
We can express $g_i$, by using some permutations of rows if necessary, as follows
\[
g_i=(\underbrace{1\cdots 1}_{\alpha}~0\cdots 0)^T.
\]
In the following we only deal with the case $\alpha<k$, but the modification to the case $\alpha=k$ is trivial.
Now we have the following two cases\vspace{0.2cm}\\
case I: there exists $l>\alpha$ such that $g_{lj}=1$,\\
case I\!I: not case I, i.e., $g_{lj}=0$ for all $l>\alpha$. \vspace{0.2cm}\\
In the case I, let us fix $(m_1\cdots m_\alpha)$ in a message $m=(m_1\cdots m_k)$ with $m_1+\cdots+m_{\alpha}=1$, 
which corresponds to $x_i=1$, and consider the numbers of codewords with $x_j=0$ and $x_j=1$ 
from the remaining message bits 
$(m_{\alpha+1} \cdots m_k)$. 
Then, from the assumption, there exists at least one non zero element in $g_{\alpha+1,j},\cdots,g_{k,j}$.
Hence, the same argument in Lemma \ref{lemma:gi} shows that the numbers of codewords with $x_j=0$ and 
$x_j=1$ under a fixed $(m_1,\cdots,m_{\alpha})$ are the same. 
By considering all the possibilities of $(m_1\cdots m_\alpha)$ with $m_1+\cdots+m_\alpha=1$, it gives 
$\sharp \mathcal{C}(\substack{x_i=1\\x_j=0})=\sharp \mathcal{C}(\substack{x_i=1\\x_j=1})=2^{k-2}.$

Next, let us consider the case I\!I. Again by using a permutation if necessary, we have the following expressions
\[
g_i=(\underbrace{1\cdots 1}_{\alpha_1}~\underbrace{1\cdots 1}_{\alpha_2}~0\cdots 0)^T,~~~
g_j=(\underbrace{1\cdots 1}_{\alpha_1}~0\cdots 0)^T,
\]
where $\alpha=\alpha_1+\alpha_2$ and $\alpha_2>0$ because of $g_i\neq g_j$. Then we have
\begin{eqnarray*}
\sharp \mathcal{C}(\substack{x_i=1\\x_j=0})&=&
\sharp\left\{m_1+\cdots+m_{\alpha_1+\alpha_2}=1~{\rm and}~
m_1+\cdots+m_{\alpha_1}=0
\right\}\times 2^{k-\alpha}\\
&=&\sharp\left\{
m_{\alpha_1+1}+\cdots+m_{\alpha_1+\alpha_2}=1~{\rm and}~
m_1+\cdots+m_{\alpha_1}=0
\right\}\times 2^{k-\alpha},\\
\sharp \mathcal{C}(\substack{x_i=1\\x_j=1})&=&
\sharp\left\{m_1+\cdots+m_{\alpha_1+\alpha_2}=1~{\rm and}~
m_1+\cdots+m_{\alpha_1}=1
\right\}\times 2^{k-\alpha}\\
&=&\sharp\left\{
m_{\alpha_1+1}+\cdots+m_{\alpha_1+\alpha_2}=0~{\rm and}~
m_1+\cdots+m_{\alpha_1}=1
\right\}\times 2^{k-\alpha}.
\end{eqnarray*}
However, the same argument in Lemma \ref{lemma:gi} implies 
$\sharp \mathcal{C}(\substack{x_i=1\\x_j=0})=\sharp \mathcal{C}(\substack{x_i=1\\x_j=1})=2^{k-2}.$
}\qed
\end{prf}\vspace{0.1cm}
{\underline {\it Proof of Proposition \ref{pp:jacobi}.}}~~~
The $(i,j)$ element in the Jacobi matrix $J$ is given by
\[
\frac{\partial f_i}{\partial u_j}(p)=\left\{
I^j_i(p)H(p)-I_i(p)H^j(p)
\right\}/H^2(p).
\]
From Lemma \ref{lemma:gi}, $H^j(p)=0$ because
\[
H^j(p)=\left(\frac{1}{2}\right)^{n-1}\!\!\!\!\!\!\times 2^{k-1}-\left(\frac{1}{2}\right)^{n-1}\!\!\!\!\!\!\times 2^{k-1}=0,
\]
where the first term comes from the codewords with $x_j=1$ and the second term comes from the codewords with $x_j=0$. 
It is also easy to observe that $H(p)=(1/2)^n\times 2^k=2^{-n+k}$, and $I^i_i(p)=(1/2)^{n-1}\times 2^{k-1}=2^{-n+k}$. 
Therefore the diagonal elements are $J_{ii}=1,i=1,\cdots,n$.

Next, let us consider the case $i\neq j$. In this case, if we have $g_i=g_j$, then, from Lemma \ref{lemma:gij}, $I^j_i(p)=(1/2)^{n-1}\times 2^{k-1}=2^{-n+k}$. On the other hand, if we have $g_i\neq g_j$, then 
$I^j_i(p)=(1/2)^{n-1}\times2^{k-2}-(1/2)^{n-1}\times2^{k-2}=0$. This concludes the proof of Proposition \ref{pp:jacobi}.
\qed\vspace{0.1cm}

Two corollaries follow from Proposition \ref{pp:jacobi} which characterize the eigenvalues and 
the eigenvectors of the Jacobi matrix $J$, and it clearly determines the stability and the stable/unstable subspaces of $p$. 
To this end, let us denote by $\mathcal{G}_J$ a graph whose adjacent matrix is $J$. 
Namely, the nodes of $\mathcal{G}_J$ are $1,\cdots,n$, and an undirected edge $(i,j)$ appears in $\mathcal{G}_J$
if and only if $J_{ij}=1$. 
\begin{corollary}\label{corollary:eigenvalue}
Suppose the graph $\mathcal{G}_J$ is decomposed into $l$ connected components
\[
\mathcal{G}_J=A_1\cup \cdots \cup A_l,~~~A_i\cap A_j=\emptyset~{\rm if}~i\neq j.
\]
Let $n_i$  be the number of nodes in the component $A_i$, $i=1,\cdots,l$. 
Then all the eigenvalues of $J$ are given by 
\[
n_1, n_2,\cdots, n_l, 0,
\]
where the eigenvalues $n_1,\cdots,n_l$ are simple and the $0$ eigenvalue has $(n_1+\cdots+n_l-l)$ multiplicity.
\end{corollary}
\begin{prf}{\rm 
From Proposition \ref{pp:jacobi}, any two nodes in a same connected component have an edge between them.
Hence, it is possible to transform $J$ into the following block diagonal matrix
\begin{equation}
E^{-1}JE=\left(\begin{array}{cccc}
B_1 & 0 & \cdots & 0\\
0 & B_2 & \cdots &0\\
\vdots & \vdots & \ddots & \vdots \\
0 & 0 & \cdots & B_l
\end{array}
\right),   \label{eq:elementarymatrix}
\end{equation}
where $E$ is determined by compositions of column switching elementary matrices, and 
$B_i$ is an $n_i\times n_i$ matrix all of whose elements are 1. The statement of the corollary follows immediately.
}\qed
\end{prf}

From now on, we treat a Jacobi matrix of the block diagonal form (\ref{eq:elementarymatrix}). 
Obviously, it gives no restriction since, if necessary, we can appropriately permute columns 
of the original generator matrix in advance. 
Let us denote the set of eigenvalues of $J$ derived in Corollary \ref{corollary:eigenvalue} by
\begin{equation}
n_1,0,\cdots,0,\cdots,n_i,0,\cdots,0,\cdots,n_l,0,\cdots,0,  \label{eq:eigenvalue}
\end{equation}
where the successive $0,\cdots, 0$ after each $n_i$ has $(n_i-1)$ elements. 
In case of $n_i=1$, we ignore the successive $0,\cdots,0$ for $n_i$ (i.e., $n_i, n_{i+1}, 0, \cdots$).

\begin{corollary}\label{corollary:eigenvector}
The Jacobi matrix $J$ is diagonalizable. Furthermore, the corresponding eigenvectors
\[
p^{(1)}_1,\cdots,p^{(1)}_{n_1},\cdots,p^{(i)}_1\cdots,p^{(i)}_{n_i},\cdots,p^{(l)}_1,\cdots,p^{(l)}_{n_l},
\]
for {\rm (\ref{eq:eigenvalue})} under this ordering are given by the following
\[
p^{(i)}_1=\left(\begin{array}{c}
\bf{0}\\
q^{(i)}_1\\
\bf{0}
\end{array}
\right),
\cdots,
p^{(i)}_{n_i}=\left(\begin{array}{c}
\bf{0}\\
q^{(i)}_{n_i}\\
\bf{0}
\end{array}
\right),
\]
where 
\[
q^{(i)}_1=\left(\begin{array}{r}
1\\
\vdots\vspace{-0.2cm}\\
\vdots\vspace{-0.2cm}\\
\vdots\\
1
\end{array}
\right),~
q^{(i)}_2=\left(\begin{array}{r}
1\\
-1\\
0\\
\vdots\\
0
\end{array}
\right),~
q^{(i)}_3=\left(\begin{array}{r}
0\\
1\\
-1\\
\vdots\\
0
\end{array}
\right),
\cdots,
q^{(i)}_{n_i}=\left(\begin{array}{r}
0\\
\vdots\\
0\\
1\\
-1
\end{array}
\right).
\]
Here the vectors $q^{(i)}_1\cdots q^{(i)}_{n_i}$ have $n_i$ elements and the first element starts at 
the $(n_1+\cdots+n_{i-1}+1)$-th row.
The bold type ${\bf 0}$ expresses that the remaining elements in $p^{(i)}_1,\cdots,p^{(i)}_{n_i}$ are $0$. 
In case of $n_i=1$, we only have $p^{(i)}_1$ {\rm (}with $q^{(i)}_1=(1)${\rm )}.
\end{corollary}
\begin{prf}{\rm
It is obvious from Corollary \ref{corollary:eigenvalue}.
}\qed
\end{prf}

We finally mention a relationship to a continuous gradient dynamical system. 
Let us denote by $F^{(i)}_{x}(u)$ a polynomial obtained by removing $\rho_i(u_i)$ from a codeword 
polynomial $F_x(u)$ of $x\in\mathcal{C}$. By using this notation, the map (\ref{eq:map}) can be also described as
\[
u'_i=f_i(u)=\frac{\sum_{\substack{x\in\mathcal{C}\\x_i=1}}F_x(u)}{\sum_{\substack{x\in\mathcal{C}\\x_i=0}}F_x(u)+\sum_{\substack{x\in\mathcal{C}\\x_i=1}}F_x(u)}=
\frac{\sum_{\substack{x\in\mathcal{C}\\x_i=1}}u_iF^{(i)}_x(u)}{\sum_{\substack{x\in\mathcal{C}\\x_i=0}}(1-u_i)F^{(i)}_x(u)+\sum_{\substack{x\in\mathcal{C}\\x_i=1}}u_iF^{(i)}_x(u)}.
\]
Then it follows that 
\begin{eqnarray*}
u'_i-u_i&=&\frac{1}{H(u)}\left(
\sum_{\substack{x\in\mathcal{C}\\x_i=1}}u_iF^{(i)}_x(u)
-\sum_{\substack{x\in\mathcal{C}\\x_i=0}}u_i(1-u_i)F^{(i)}_x(u)
-\sum_{\substack{x\in\mathcal{C}\\x_i=1}}u^2_iF^{(i)}_x(u)
\right)\\
&=&\frac{u_i(1-u_i)}{H(u)}\left(
\sum_{\substack{x\in\mathcal{C}\\x_i=1}}F^{(i)}_x(u)-\sum_{\substack{x\in\mathcal{C}\\x_i=0}}F^{(i)}_x(u)
\right)\\
&=&u_i(1-u_i)\frac{\partial}{\partial u_i}\left(\log H(u)\right).
\end{eqnarray*}
This proves the following proposition.
\begin{pp}\label{pp:gradient}
The rational map {\rm (\ref{eq:map})} maps a point $u\in {\rm Int}(I^n)$ to the direction of ${\rm grad}(\log H(u))$ 
with a contraction rate $u_i(1-u_i)$ for each element $i=1,\cdots,n$. Especially, $u\in {\rm Int}(I^n)$ is a fixed point
of {\rm (\ref{eq:map})} if and only if it is a fixed point of the continuous gradient dynamical system 
$du/dt={\rm grad}(\log H(u))$.
\end{pp}

\section{Generator Matrix and Taylor Expansion}\label{sec:taylor}
Next, we study a relationship between higher order terms in (\ref{eq:taylor}) and a generator matrix $G$, 
and prove Theorem \ref{thm:thm}. 
For this purpose, the key proposition is given as follows.
\begin{pp} \label{pp:fl0}
Let $l\geq 2$. 
If any distinct $(l+1)$ column vectors of $G$ are linearly independent, then $f^{(l)}(v)=0$.
\end{pp}

Let us denote a higher order derivative of an $i$-th element $f_{i}$ with respect to variables 
$u_{i_1}, \cdots, u_{i_l}$ by
\[
f^{i_1\cdots i_l}_{i}(p):=\frac{\partial^l f_{i}}{\partial u_{i_1}\cdots \partial u_{i_l}}(p).
\]
For the proof of Proposition \ref{pp:fl0}, we need to study higher order derivatives $H^{i_1\cdots i_l}(p)$ and 
$I^{i_1\cdots i_l}_{i}(p)$.
%
Let us at first focus on higher order derivatives $H^{i_1\cdots i_l}(p)$. We begin with the following observation, 
which characterizes the numbers of subcodes by means of column vectors of $G$. 
It is noted that we adopt the $mod~2$ arithmetic for elements in $\mathbb{F}_2$. 
\begin{lemma}\label{lemma:C01}
Suppose $i_1,\cdots,i_l\in\{1,\cdots,n\},~l\geq 2,$ are distinct indices, and let 
\begin{eqnarray*}
&&\mathcal{C}(x_{i_1}+\cdots+x_{i_l}=0)=\left\{x\in \mathcal{C}~|~x_{i_1}+\cdots+x_{i_l}=0\right\},\\
&&\mathcal{C}(x_{i_1}+\cdots+x_{i_l}=1)=\left\{x\in \mathcal{C}~|~x_{i_1}+\cdots+x_{i_l}=1\right\}
\end{eqnarray*}
be subcodes in $\mathcal{C}$.
Then the followings hold
\begin{enumerate}
\item if $g_{i_1}+\cdots+g_{i_l}\neq 0$, then $\sharp\mathcal{C}(x_{i_1}+\cdots+x_{i_l}=0)=\sharp\mathcal{C}(x_{i_1}+\cdots+x_{i_l}=1)=2^{k-1}$.
\item if $g_{i_1}+\cdots+g_{i_l}= 0$, then $\sharp\mathcal{C}(x_{i_1}+\cdots+x_{i_l}=0)=2^k$~and~
${\sharp\mathcal{C}(x_{i_1}+\cdots+x_{i_l}=1)=0}$.
\end{enumerate}
\end{lemma}
\begin{prf}{\rm 
By a suitable bit permutation, if necessary, the sum of $g_{i_1},\cdots,g_{i_l}$ can be expressed as
\[
g_{i_1}+\cdots+g_{i_l}=(\underbrace{1\cdots 1}_{\alpha}~0\cdots 0)^T.
\]
Then an original message $(m_1 \cdots m_k)$ and its codeword $(x_1 \cdots x_n)$ satisfy the following
\[
x_{i_1}+\cdots+x_{i_l}=(m_1~\cdots~m_k)\cdot(g_{i_1}+\cdots+g_{i_l}).
\]

In the case 1 ($\alpha\neq 0$), it leads to $x_{i_1}+\cdots+x_{i_l}=m_1+\cdots+m_{\alpha}$,  
so the conclusion follows from the same argument in Lemma \ref{lemma:gi}.
The case 2 is trivial from the above expression of $x_{i_1}+\cdots+x_{i_l}$.
}\qed
\end{prf}
The following lemma classifies the value $H^{i_1\cdots i_l}(p)$ based on the column vectors of $G$.
\begin{lemma}\label{lemma:H}
Let $l\geq 2$. Then $H^{i_1\cdots i_l}(p)=0$ if either of 
\begin{enumerate}
\item there exist same indices in $i_1,\cdots,i_l$
\item $g_{i_1}+\cdots+g_{i_l}\neq 0$
\end{enumerate}
is satisfied. Otherwise, that is $i_1,\cdots,i_l$ are all distinct and $g_{i_1}+\cdots+g_{i_l}=0$, 
\[
H^{i_1\cdots i_l}(p)=(-1)^l 2^{-n+k+l}.
\]
\end{lemma}
\begin{prf}{\rm 
The condition 1 immediately implies $H^{i_1\cdots i_l}(p)=0$ since the degree of each variable $u_i$ in $H(u)$ is 1. 
Hence we assume all the indices are distinct. 

Let us define the following subcodes
\begin{eqnarray*}
&&Z^{\rm odd}_{i_1\cdots i_l}:=\{x\in \mathcal{C}~|~{\rm the~number~of~0~in}~x_{i_1},\cdots,x_{i_l}~{\rm is~odd}\}\\
&&Z^{\rm even}_{i_1\cdots i_l}:=\{x\in \mathcal{C}~|~{\rm the~number~of~0~in}~x_{i_1},\cdots,x_{i_l}~{\rm is~even}\}.
\end{eqnarray*}
Then $H^{i_1\cdots i_l}(p)$ can be expressed by
\[
H^{i_1\cdots i_l}(p)=(1/2)^{n-l}(\sharp Z^{\rm even}_{i_1\cdots i_l}-\sharp Z^{\rm odd}_{i_1\cdots i_l}).
\]
Suppose $g_{i_1}+\cdots+g_{i_l}\neq 0$. Then, by Lemma \ref{lemma:C01}, we have
$\sharp \mathcal{C}(x_{i_1}+\cdots+x_{i_l}=0)=\sharp \mathcal{C}(x_{i_1}+\cdots+x_{i_l}=1)$. 
On the other hand, when $l$ is odd (or even, resp.), 
\begin{eqnarray*}
&&Z^{\rm odd}_{i_1\cdots i_l}=\mathcal{C}(x_{i_1}+\cdots+x_{i_l}=0) ~({\rm or}~ =\mathcal{C}(x_{i_1}+\cdots+x_{i_l}=1), {\rm resp.}),\\
&&Z^{\rm even}_{i_1\cdots i_l}=\mathcal{C}(x_{i_1}+\cdots+x_{i_l}=1) ~({\rm or}~=\mathcal{C}(x_{i_1}+\cdots+x_{i_l}=0), {\rm resp.}).
\end{eqnarray*}
Therefore it concludes $H^{i_1\cdots i_l}(p)=0$. 
The statement for $g_{i_1}+\cdots+g_{i_l}=0$ is similarly derived from Lemma \ref{lemma:C01}.
}\qed
\end{prf}

Next, we try to classify the value $I^{i_1\cdots i_l}_{i}(p)$. 
\begin{lemma}\label{lemma:iC01}
Suppose $i, i_1,\cdots,i_l\in\{1,\cdots,n\}, l\geq 2,$ are distinct indices and let us define two subcodes 
in $\mathcal{C}$ by
\begin{eqnarray*}
&&\mathcal{C}\left(\substack{x_{i}=1 \\ x_{i_1}+\cdots+x_{i_l}=0}\right):=\left\{x\in\mathcal{C}~|~x_{i}=1,~x_{i_1}+\cdots+x_{i_l}=0\right\},\\
&&\mathcal{C}\left(\substack{x_{i}=1 \\ x_{i_1}+\cdots+x_{i_l}=1}\right):=\left\{x\in\mathcal{C}~|~x_{i}=1,~x_{i_1}+\cdots+x_{i_l}=1\right\}.
\end{eqnarray*}
Then the following classification holds
\begin{enumerate}
\item if $g_{i_1}+\cdots+g_{i_l}=0$, then $\sharp\mathcal{C}\left(\substack{x_{i}=1 \\ x_{i_1}+\cdots+x_{i_l}=0}\right)=2^{k-1}$,  $\sharp\mathcal{C}\left(\substack{x_{i}=1 \\ x_{i_1}+\cdots+x_{i_l}=1}\right)=0$.
\item if $0\neq g_{i_1}+\cdots+g_{i_l}\neq g_{i}$, then $\sharp\mathcal{C}\left(\substack{x_{i}=1 \\ x_{i_1}+\cdots+x_{i_l}=0}\right)=\sharp\mathcal{C}\left(\substack{x_{i}=1 \\ x_{i_1}+\cdots+x_{i_l}=1}\right)=2^{k-2}$.
\item if $0\neq g_{i_1}+\cdots+g_{i_l}=g_{i}$, then $\sharp\mathcal{C}\left(\substack{x_{i}=1 \\ x_{i_1}+\cdots+x_{i_l}=0}\right)=0$, $\sharp\mathcal{C}\left(\substack{x_{i}=1 \\ x_{i_1}+\cdots+x_{i_l}=1}\right)=2^{k-1}$.
\end{enumerate}
\end{lemma}
\begin{prf}{\rm
Since we have
\[
x_{i_1}+\cdots+x_{i_l}=(m_1~\cdots~m_k)\cdot(g_{i_1}+\cdots+g_{i_l}),
\]
the case 1 and 3 are trivial. So, let us assume $0\neq g_{i_1}+\cdots+g_{i_l}\neq g_{i}$. The proof is similar to that of Lemma 
\ref{lemma:gij}. By using a suitable permutation, let us express $g_{i}$ as follows
\[
g_{i}=(\underbrace{1\cdots 1}_{\alpha}~0\cdots 0)^T.
\]
Here we only deal with the case $\alpha<k$ again, since the modification for $\alpha = k$ follows immediately from the 
following case I\!I. We have two situations\vspace{0.2cm}\\
case I: there exists $\beta>\alpha$ such that $g_{\beta,i_1}+\cdots+g_{\beta,i_l}=1$.\\
case I\!I: not case I, i.e., $g_{\beta,i_1}+\cdots+g_{\beta,i_l}=0$ for all $\beta>\alpha$.\vspace{0.2cm}\\
In case I, let us fix $(m_1~\cdots~m_\alpha)$ with $m_1+\cdots+m_{\alpha}=1$, which corresponds to $x_{i}=1$, 
and consider the numbers of codewords with $x_{i_1}+\cdots+x_{i_l}=0$ or $=1$ for the remaining message bits
$(m_{\alpha+1}~\cdots~m_k)$. From the assumption, the $\beta$-th element of the vector $g_{i_1}+\cdots+g_{i_l}$
is 1, and, by applying Lemma \ref{lemma:C01} to the subvector from $(\alpha+1)$-th to $k$-th elements, 
the numbers of codewords with $x_{i_1}+\cdots+x_{i_l}=0$ or $=1$ are the same for each $(m_1~\cdots~m_\alpha)$.
Hence we have $\sharp\mathcal{C}\left(\substack{x_{i}=1 \\ x_{i_1}+\cdots+x_{i_l}=0}\right)=\sharp\mathcal{C}\left(\substack{x_{i}=1 \\ x_{i_1}+\cdots+x_{i_l}=1}\right)=2^{k-2}$.
The proof for case I\!I is almost parallel to that of Lemma \ref{lemma:gij}, so we omit it.
}\qed
\end{prf}

Let us introduce the following notations, which are similar to those in the proof of Lemma \ref{lemma:H},
\begin{eqnarray*}
&&\hspace{-0.5cm}Z^{\rm odd}_{i_1\cdots i_l}(x_{i}=1)=\{x\in \mathcal{C}~|~{\rm the~number~of~0~in}~x_{i_1},\cdots,x_{i_l}~{\rm is~odd~and}~x_{i}=1\},\\
&&\hspace{-0.5cm}Z^{\rm even}_{i_1\cdots i_l}(x_{i}=1)=\{x\in \mathcal{C}~|~{\rm the~number~of~0~in}~x_{i_1},\cdots,x_{i_l}~{\rm is~even~and}~x_{i}=1\}.
\end{eqnarray*}
Then, for odd $l$ (or even $l$, resp.),  we have 
\begin{eqnarray*}
&&Z^{\rm odd}_{i_1\cdots i_l}(x_{i}=1)=\mathcal{C}\left(\substack{x_{i}=1 \\ x_{i_1}+\cdots+x_{i_l}=0}\right)~
({\rm or~}=\mathcal{C}\left(\substack{x_{i}=1 \\ x_{i_1}+\cdots+x_{i_l}=1}\right),~{\rm resp.}),\\
&&Z^{\rm even}_{i_1\cdots i_l}(x_{i}=1)=\mathcal{C}\left(\substack{x_{i}=1 \\ x_{i_1}+\cdots+x_{i_l}=1}\right)~
({\rm or~}=\mathcal{C}\left(\substack{x_{i}=1 \\ x_{i_1}+\cdots+x_{i_l}=0}\right),~{\rm resp.}).
\end{eqnarray*}

%
\begin{lemma}\label{lemma:Ip}
Let $l\geq 2$. Then $I^{i_1\cdots i_l}_{i}(p)$ is classified as
\[
I^{i_1\cdots i_l}_{i}(p)=\left\{
\begin{array}{ll}
0, & {\rm If~C0}~or~{\rm C1}\\
(-1)^{l+1}2^{-n+k+l-1}, & {\rm If~C2}\\
(-1)^{l}2^{-n+k+l-1}, & {\rm If~C3}
\end{array}
\right.,
\]
where each condition is given by\vspace{0.3cm}\\
{\rm C0:} there exist same indices in $i_1,\cdots,i_l$\vspace{0.1cm}\\
{\rm C1:} $\overline{{\rm C0}}$ and $0\neq g_{i_1}+\cdots+g_{i_l}\neq g_{i}$~~(here $\overline{{\rm C0}}$ means 
``NOT {\rm C0}")\vspace{0.1cm}\\
{\rm C2:} $\overline{{\rm C0}}$ and $g_{i_1}+\cdots+g_{i_l}=g_{i}$\vspace{0.1cm}\\
{\rm C3:} $\overline{{\rm C0}}$ and $g_{i_1}+\cdots+g_{i_l}=0$
\end{lemma}
\begin{prf}{\rm
The condition C0 immediately implies $I^{i_1\cdots i_l}_{i}(p)=0$. Hence we assume all the indices are distinct.
The remaining proof follows directly from Lemma \ref{lemma:iC01} for each case.
First of all, let us study the case $i\notin\{i_1,\cdots,i_l\}$. By using the notations introduced before the lemma,
we have
\[
I^{i_1\cdots i_l}_{i}(p)=(1/2)^{n-l}(\sharp Z^{\rm even}_{i_1\cdots i_l}(x_{i}=1)-\sharp Z^{\rm odd}_{i_1\cdots i_l}(x_{i}=1)).
\]
Therefore the condition C1 implies $I^{i_1\cdots i_l}_{i}(p)=0$ by Lemma \ref{lemma:iC01}.

On the other hand, if we assume the condition C2, then $\sharp\mathcal{C}\left(\substack{x_{i}=1 \\ x_{i_1}+\cdots+x_{i_l}=0}\right)=0$ and $\sharp\mathcal{C}\left(\substack{x_{i}=1 \\ x_{i_1}+\cdots+x_{i_l}=1}\right)=2^{k-1}$ from 
Lemma \ref{lemma:iC01}. 
If $l$ is even, $I^{i_1\cdots i_l}_{i}(p)=-(1/2)^{n-l}\times 2^{k-1}=-2^{-n+k+l-1}$. 
Similarly, if $l$ is odd, $I^{i_1\cdots i_l}_{i}(p)=2^{-n+k+l-1}$. 
Hence, we have $I^{i_1\cdots i_l}_{i}(p)=(-1)^{l+1}2^{-n+k+l-1}$ for the condition C2. 

For the condition C3, the role of $\sharp\mathcal{C}\left(\substack{x_{i}=1 \\ x_{i_1}+\cdots+x_{i_l}=0}\right)$ 
and $\sharp\mathcal{C}\left(\substack{x_{i}=1 \\ x_{i_1}+\cdots+x_{i_l}=1}\right)$ changes each other 
from Lemma \ref{lemma:iC01}, so it just leads to the opposite sign in $I^{i_1\cdots i_l}_{i}(p)$ to that for the condition C2.

Next, let us study the case $i\in\{i_1,\cdots,i_l\}$. Without loss of generality, let us suppose $i=i_l$. Then we have
\[
I^{i_1\cdots i_l}_{i}(p)=(1/2)^{n-l}(\sharp Z^{\rm even}_{i_1\cdots i_{l-1}}(x_{i}=1)-\sharp Z^{\rm odd}_{i_1\cdots i_{l-1}}(x_{i}=1)).
\]
Therefore the condition C1 implies $I^{i_1\cdots i_l}_{i}(p)=0$ by Lemma \ref{lemma:iC01} (or Lemma \ref{lemma:gij} for $l=2$). 

For the condition C2 ($l$ must be $>2$), 
we have $\sharp\mathcal{C}\left(\substack{x_{i}=1 \\ x_{i_1}+\cdots+x_{i_{l-1}}=0}\right)=2^{k-1}$ and 
$\sharp\mathcal{C}\left(\substack{x_{i}=1 \\ x_{i_1}+\cdots+x_{i_{l-1}}=1}\right)=0$ from Lemma \ref{lemma:iC01}. 
Hence, by the same calculation as that for $i\notin\{i_1,\cdots,i_l\}$, we have $I^{i_1\cdots i_l}_{i}(p)=(-1)^{l+1}2^{-n+k+l-1}$.

Finally, let us consider the condition C3. For $l>2$, we have $\sharp\mathcal{C}\left(\substack{x_{i}=1 \\ x_{i_1}+\cdots+x_{i_{l-1}}=0}\right)=0$ and $\sharp\mathcal{C}\left(\substack{x_{i}=1 \\ x_{i_1}+\cdots+x_{i_{l-1}}=1}\right)=2^{k-1}$
by Lemma \ref{lemma:iC01} again. So, $I^{i_1\cdots i_l}_{i}(p)=(-1)^{l}2^{-n+k+l-1}$. 
For $l=2$, we can not use Lemma \ref{lemma:iC01} because of $l-1=1<2$. 
However, $x_{i}=x_{i_1}$ holds from the assumption $g_{i}=g_{i_1}$. Hence a direct calculation shows 
$I^{i_1i_2}_{i}(p)=(1/2)^{n-2}\times 2^{k-1}=2^{-n+k+1}$, which is the formula for $l=2$. 
}\qed
\end{prf}

Before proving Proposition \ref{pp:fl0}, let us show the following two lemmas. The proofs of them are easy 
application of induction. 
\begin{lemma}\label{lemma:1/A}
The derivative $(1/H)^{i_1\cdots i_l}$ is given by
\[
\left(\frac{1}{H}\right)^{i_1\cdots i_l}\!\!\!\!\!\!=(-1)^l\frac{l!}{H^{l+1}}(H^{i_1}\cdots H^{i_l})+\cdots
+(-1)^k\frac{k!}{H^{k+1}}\sum_{C(l,k)}(H^{{\bf r}^1}\cdots H^{{\bf r}^k})+\cdots
-\frac{1}{H^2}H^{i_1\cdots i_l},
\]
where the summation for the $k$-th term {\rm (}$1<k<l${\rm )} is taken on all combinations $C(l,k)$ for 
dividing $i_1,\cdots,i_l$ into $k$ groups. 
Here ${\bf r}^1,\cdots,{\bf r}^k$ represent a decomposition of $\{i_1,\cdots,i_l\}${\rm ;}
\[
\{i_1,\cdots,i_l\}=\cup_{i=1}^{k}{\bf r}^i,~~~{\bf r}^i\neq \emptyset,~~~{\bf r}^i\cap {\bf r}^j=\emptyset~(i\neq j).
\]
\end{lemma}
\begin{lemma}\label{lemma:fl}
Let $l\geq 2$. Then the derivative $f^{i_1\cdots i_l}_{i}$ is given by
\[
f^{i_1\cdots i_l}_{i}=I^{i_1\cdots i_l}_{i}\left(\frac{1}{H}\right)
+\cdots+
\sum^{_{l}C_k}_{j=1}I_{i}^{{\bf p}^k_j}\left(\frac{1}{H}\right)^{{\bf q}^k_j}
+\cdots+
I_{i}\left(\frac{1}{H}\right)^{i_1\cdots i_l},
\]
where ${\bf p}^k_j$ and ${\bf q}^k_j$ are a decomposition of $\{i_1,\cdots,i_l\}${\rm ;}
\[
{\bf p}^k_j\cap {\bf q}^k_j=\emptyset,~~~{\bf p}^k_j\cup {\bf q}^k_j=\{i_1,\cdots,i_l\},~~~\sharp {\bf q}^k_j=k,
\]
and the summations are taken on all the combinations of the decompositions.
\end{lemma}
Now we prove Proposition \ref{pp:fl0}.\vspace{0.1cm}\\
{\underline {\it Proof of Proposition \ref{pp:fl0}.}}~~~
Let us assume that any distinct $(l+1)$ column vectors of $G$ are linearly independent. Then, from Lemma \ref{lemma:H}, 
\ref{lemma:1/A}, and \ref{lemma:fl}, we have 
\[
f^{i_1\cdots i_l}_{i}(p)=I^{i_1\cdots i_l}_{i}(p)\left(\frac{1}{H(p)}\right).
\]
On the other hand, from Lemma \ref{lemma:Ip}, $I^{i_1\cdots i_l}_{i}(p)$ is 0. The proof is completed.
\qed\vspace{0.3cm}\\
Finally, we are in the position to prove Theorem \ref{thm:thm}.\vspace{0.1cm}\\
{\underline {\it Proof of Theorem \ref{thm:thm}.}}~~~
The formula for the case $l=2$ is given by Proposition \ref{pp:jacobi}. 
Let us assume $l\geq 3$. 
From Proposition \ref{pp:fl0}, all nonlinear terms with orders less than $l$ are zero, 
so we only study the $l$-th nonlinear terms.
From the assumption and similar argument in the proof of Proposition \ref{pp:fl0}, 
the derivative $f^{i_1\cdots i_l}_{i}(p)$ is given by
\[
f^{i_1\cdots i_l}_{i}(p)=I^{i_1\cdots i_l}_{i}(p)\left(\frac{1}{H(p)}\right).
\]

Then, the classification in Lemma \ref{lemma:Ip} shows that 
$f^{i_1\cdots i_l}_{i}(p)=0$ in C0 and C1. Moreover, the condition C3 does not occur from the assumption.
Hence, only nonzero terms are derived from the condition C2.
It should be noted that the set of indices $(i_1,\cdots,i_l)$ satisfying the condition C2 with $i$ is 
exactly the same as $\Theta^{(l)}_{i}$.
Hence, Lemma \ref{lemma:Ip} results in 
\begin{eqnarray*}
f^{i_1\cdots i_l}_{i}(p)&=&I^{i_1\cdots i_l}_{i}(p)\left(\frac{1}{H(p)}\right)\\
&=&(-1)^{l+1}2^{-n+k+l-1}\frac{1}{(1/2)^n\times 2^k}\\
&=&(-2)^{l-1},
\end{eqnarray*}
for $(i_1,\cdots,i_l)\in\Theta^{(l)}_{i}$.
Finally, we have the following 
\begin{eqnarray*}
f^{(l)}_{i}(v)&=&\sum_{\substack{m_1+\cdots+m_n= l\\m_k\geq 0}}
\frac{v_1^{m_1}\cdots v_n^{m_n}}{m_1!\cdots m_n!}
\left(\frac{\partial^l f_i}{\partial u_1^{m_1}\cdots\partial u_n^{m_n}}\right)(p)\\
&=&\sum_{(i_1,\cdots,i_l)\in \Theta^{(l)}_i}(-2)^{l-1}
v_{i_1}\cdots v_{i_l}.
\end{eqnarray*}
The identity in (\ref{eq:hosum}) is derived by just substituting $v_i=u_i-1/2$.
It completes the proof.
\qed\vspace{0.3cm}
\section{Relationship to Algebraic Geometry Codes}\label{sec:ag}
The purpose of this section is to derive Corollary \ref{corollary:ag}.
This corollary shows that usual techniques in algebraic geometry codes can be applied to control,  
not only the minimum distance of the code, but also the approximate ML decoding. 
For the definitions of basic tools in algebraic geometry such as
genus, divisor, Riemann-Roch space, differential, and residue, 
we refer to \cite{fulton}.
We also request basic knowledge of algebraic geometry codes in this section 
(e.g., see \cite{vanlint}, \cite{stich} and \cite{tvn}). 

Let $\mathbb{F}_q$ be a finite field with $q>1$ elements.
Let us first recall two classes of algebraic geometry codes called geometric Reed-Solomon codes and geometric Goppa codes. 
Let $\mathcal{X}$ be an absolutely irreducible nonsingular projective curve over $\mathbb{F}_q$. 
For rational points $P_1,\cdots,P_{\tilde n}$ on $\mathcal{X}$, we define a divisor on $\mathcal{X}$ by 
$D=P_1+\cdots +P_{\tilde n}$.
Moreover, let $D'$ be another divisor whose support is disjoint to $D$. We assume that $D'$ satisfies the following condition
\[
2g-2<{\rm deg}(D')<\tilde{n}
\] 
for the sake of simplicity. Here $g$ is the genus of $\mathcal{X}$.
Geometric Reed-Solomon codes are characterized by the Riemann-Roch space associated to $D'$ 
\[
\mathcal{L}(D'):=\left\{
\phi\in\mathbb{F}^*_q(\mathcal{X})~|~(\phi)+D'\geq 0
\right\}\cup\{0\},
\]
where $\mathbb{F}^*_q(\mathcal{X})$ is the set of nonzero elements of the function field $\mathbb{F}_q(\mathcal{X})$, 
and $(\phi)$ is the principal divisor of the rational function $\phi$.
The fundamental fact that the Riemann-Roch space $\mathcal{L}(D')$ is a finite dimensional vector space leads to the 
following definition.
\begin{df}{\rm 
{\it The geometric Reed-Solomon code} $\mathcal{C}(D,D')$ of length $\tilde{n}$ over $\mathbb{F}_q$ is defined by 
the image of the linear map $\alpha:\mathcal{L}(D')\rightarrow\mathbb{F}^{\tilde n}_q$
given by $\alpha(\phi)=(\phi(P_1),\cdots,\phi(P_{\tilde n}))$. 
}
\end{df}
On the other hand, geometric Goppa codes are defined via differentials and their residues. 
Let us denote the set of differentials on $\mathcal{X}$ by $\Omega(\mathcal{X})$, and define for each divisor $D$
\[
\Omega(D):=\{\omega\in\Omega(\mathcal{X})~|~(\omega)-D\geq 0\},
\]
where $(\omega)$ is the divisor of the differential $\omega$. 
\begin{df}{\rm 
{\it The geometric Goppa code} $\mathcal{C}^*(D,D')$ of length ${\tilde n}$ over $\mathbb{F}_q$ is defined by
the image of the linear map $\alpha^*:\Omega(D'-D)\rightarrow \mathbb{F}^{\tilde n}_q$ given by
$\alpha^*(\omega)=({\rm Res}_{P_1}(\omega),\cdots,{\rm Res}_{P_{\tilde n}}(\omega))$, where 
${\rm Res}_{P}(\omega)$ expresses the residue of $\omega$ at $P$.
}
\end{df}

The following propositions are an easy consequence of the Riemann-Roch theorem.
\begin{pp}{\rm (e.g., \cite{vanlint}, \cite{stich}, \cite{tvn})}\label{pp:grs}
\begin{enumerate}
\item The dimension of the geometric Reed-Solomon code $\mathcal{C}(D,D')$ is $k={\rm deg}(D')-g+1$ and 
the minimum distance satisfies $d\geq \tilde{n}-{\rm deg}(D')$.
\item The dimension of the geometric Goppa code $\mathcal{C}^*(D,D')$ is $k=\tilde{n}-{\rm deg}(D')+g-1$ and 
the minimum distance satisfies $d\geq {\rm deg}(D')-2g+2$.
\item The codes $\mathcal{C}(D,D')$ and $\mathcal{C}^*(D,D')$ are dual codes.
\end{enumerate}
\end{pp}
It should be noted that the minimum distances of these two codes are controlled by the genus of $\mathcal{X}$ and 
the choice of the divisor $D'$, and,  as a result, they induce appropriate linear independence on their parity check matrices.

For an application of Theorem \ref{thm:thm}, we need to derive expanded codes over $\mathbb{F}_2$ from 
geometric Reed-Solomon and geometric Goppa codes over $\mathbb{F}_{q}, q=2^s$. 
Let $\mathcal{C}_{q}$ be a code over $\mathbb{F}_q$ with length $\tilde{n}$, and $e_1,\cdots,e_s\in\mathbb{F}_q$ 
be a basis of $\mathbb{F}_2$-vector space $\mathbb{F}_q$. 
This basis naturally induces the map $\mathbb{F}^{\tilde n}_q\ni x\mapsto \hat{x}\in \mathbb{F}^{{\tilde n}s}_2$ by 
expressing each element $x_i$ in $x=(x_1,\cdots,x_{\tilde n})$ as coefficients of $\mathbb{F}_2$-vector space.
Then the expanded code of $\mathcal{C}_{q}$ over $\mathbb{F}_2$ is defined by
$
\mathcal{C}_2:=\{\hat{x}~|~x\in \mathcal{C}_q\}.
$
A relationship between $\mathcal{C}_q$ and $\mathcal{C}_2$ is given by the following proposition. 
\begin{pp}
If a code $\mathcal{C}_q$ has parameters $[\tilde{n},k,d]$, where $\tilde{n}$ is the code length, $k$ is the dimension, and 
$d$ is the minimum distance, then its expanded code $\mathcal{C}_2$ has the parameters $[\tilde{n}s, ks, d'\geq d]$.
\end{pp}

Now we apply Theorem \ref{thm:thm} to geometric Reed-Solomon/Goppa codes. 
Let $n=\tilde{n}s$, and let $\mathcal{C}_2(D,D')$ and $\mathcal{C}^*_2(D,D')$ be the expanded codes 
over $\mathbb{F}_2$ of
a geometric Reed-Solomon code  $\mathcal{C}(D,D')$ and a geometric Goppa code $\mathcal{C}^*(D,D')$ 
over $\mathbb{F}_q$. Then, we have the following corollary.
\begin{corollary}\label{corollary:ag}
The expanded geometric Reed-Solomon code $\mathcal{C}_2(D,D')$ has the minimum distance 
$d\geq {\tilde n}-{\rm deg}(D')$. 
Furthermore, there exists an $l$-th order approximate ML decoding with $l\geq{\rm deg}(D')-2g+1$ such that
$\tilde{f}_{i}(u)=u_i+f^{(l)}_i(v)$.
\end{corollary}
\begin{prf}{\rm 
The first statement is the property of a geometric Reed-Solomon code and its expansion. 
The second statement follows from Theorem \ref{thm:thm} and the duality of $\mathcal{C}(D,D')$ and $\mathcal{C}^*(D,D')$.
}\qed
\end{prf}
%
%
%
%
%
\vspace{0.3cm}
{\bf Example: Hermitian Code} (e.g., \cite{vanlint}, \cite{stich}, \cite{tvn})\\
%
Let $q=r^2$ be a power of 2. 
The Hermitian curve $\mathcal{H}$ is given by the homogeneous equation $X^{r+1}+Y^{r+1}+Z^{r+1}=0$ and 
its genus is $g=r(r-1)/2$, because there are no singular points. 
It is known that the number of rational points over $\mathbb{F}_q$ is $r^3+1$.

Let us fix $r=2$ as an example. Then, the following is the list of the rational points on $\mathcal{H}$:
\begin{eqnarray*}
&&P_1=(1,0,\bar{\omega}),~~P_2(1,0,\omega),~~P_3(1,0,1),~~P_4=(1,\bar{\omega},0),~~P_5(1,\omega,0),\\
&&P_6(1,1,0),~~P_7=(0,\bar{\omega},1),~~P_8(0,\omega,1),~~Q=(0,1,1),
\end{eqnarray*}
where $\omega$ is a primitive element of $\mathbb{F}_4$ and $\bar{\omega}=1+\omega$ 
(i.e., $\mathbb{F}_4=\{0,1,\omega,\bar{\omega}\}$).
Let us suppose $D=P_1+\cdots+P_8$ (hence the code length is ${\tilde n}=8$), and $D'=mQ$, $2g-2<m<{\tilde n}$.
A basis of the Riemann-Roch space $\mathcal{L}(D')$ with $m=4$ is given by 
\[
\mathcal{L}(4Q)={\rm Span}\left\{
1, \frac{X}{Y+Z},\frac{Y}{Y+Z},\frac{X^2}{(Y+Z)^2}
\right\}.
\]
Then, we can explicitly show a generator matrix of the expanded geometric Reed-Solomon code 
$\mathcal{C}(D,D')$ over $\mathbb{F}_2$ as follows
\begin{equation}
G=\left(
\begin{array}{cccccccccccccccc}
1 & 0 & 0 & 0 & 0 & 0 & 0 & 0 & 1 & 0 & 1 & 0 & 0 & 1 & 1 & 1\\
0 & 1 & 0 & 0 & 0 & 0 & 0 & 0 & 0 & 1 & 0 & 1 & 1 & 1 & 1 & 0\\
0 & 0 & 1 & 0 & 0 & 0 & 0 & 0 & 1 & 0 & 0 & 0 & 1 & 0 & 1 & 0\\
0 & 0 & 0 & 1 & 0 & 0 & 0 & 0 & 0 & 1 & 0 & 0 & 0 & 1 & 0 & 1\\
0 & 0 & 0 & 0 & 1 & 0 & 0 & 0 & 0 & 0 & 1 & 0 & 1 & 0 & 1 & 0\\
0 & 0 & 0 & 0 & 0 & 1 & 0 & 0 & 0 & 0 & 0 & 1 & 0 & 1 & 0 & 1\\
0 & 0 & 0 & 0 & 0 & 0 & 1 & 0 & 1 & 0 & 1 & 0 & 1 & 1 & 0 & 1\\
0 & 0 & 0 & 0 & 0 & 0 & 0 & 1 & 0 & 1 & 0 & 1 & 1 & 0 & 1 & 1
\end{array}
\right). \label{eq:Gexample}
\end{equation}
It should be noted that this case leads to a self-dual code ${\mathcal{C}(D,D')=\mathcal{C}^*(D,D')}$.
Hence, its parity check matrix $H$ is the same as $G$, and a direct calculation proves that the minimum distance of this 
code is 4. 
It means that we have the $3$-rd order ML decoding in Corollary \ref{corollary:ag}. 
The explicit forms of $\tilde{f}_{i}(u), i=1,\cdots,16,$ are given as
\begin{eqnarray*}
&&\hspace{-0.6cm}
\tilde{f}_{1}(u)=u_1 + 4v_{7}(v_{3}v_{9} + v_{5}v_{11} + v_{13}v_{15}),~~~
\tilde{f}_{2}(u)=u_2 + 4v_{8}(v_{4}v_{10} + v_{6}v_{12} + v_{14}v_{16}),\\
&&\hspace{-0.6cm}
\tilde{f}_{3}(u)=u_3 + 4v_{9}(v_{1}v_{7} + v_{5}v_{11} + v_{13}v_{15}),~~~
\tilde{f}_{4}(u)=u_4 + 4v_{10}(v_{2}v_{8} + v_{6}v_{12} + v_{14}v_{16}),\\
&&\hspace{-0.6cm}
\tilde{f}_{5}(u)=u_5 + 4v_{11}(v_{1}v_{7} + v_{3}v_{9} + v_{13}v_{15}),~~~
\tilde{f}_{6}(u)=u_6 + 4v_{12}(v_{2}v_{8} + v_{4}v_{10} + v_{14}v_{16}),\\
&&\hspace{-0.6cm}
\tilde{f}_{7}(u)=u_7 + 4v_{1}(v_{3}v_{9} + v_{5}v_{11} + v_{13}v_{15}),~~~
\tilde{f}_{8}(u)=u_8 + 4v_{2}(v_{4}v_{10} + v_{6}v_{12} + v_{14}v_{16}),\\
&&\hspace{-0.6cm}
\tilde{f}_{9}(u)=u_9 + 4v_{3}(v_{1}v_{7} + v_{5}v_{11} + v_{13}v_{15}),~~~
\tilde{f}_{10}(u)=u_{10} + 4v_{4}(v_{2}v_{8} + v_{6}v_{12} + v_{14}v_{16}),\\
&&\hspace{-0.6cm}
\tilde{f}_{11}(u)=u_{11} + 4v_{5}(v_{1}v_{7} + v_{3}v_{9} + v_{13}v_{15}),~~
\tilde{f}_{12}(u)=u_{12} + 4v_{6}(v_{2}v_{8} +4v_{4}v_{10} + v_{14}v_{16}),\\
&&\hspace{-0.6cm}
\tilde{f}_{13}(u)=u_{13} + 4v_{15}(v_{1}v_{7} + v_{3}v_{9} + v_{5}v_{11}),~~
\tilde{f}_{14}(u)=u_{14} + 4v_{16}(v_{2}v_{8} + v_{4}v_{10} + v_{6}v_{12}),\\
&&\hspace{-0.6cm}
\tilde{f}_{15}(u)=u_{15} + 4v_{13}(v_{1}v_{7} + v_{3}v_{9} + v_{5}v_{11}),~~
\tilde{f}_{16}(u)=u_{16} + 4v_{14}(v_{2}v_{8} + v_{4}v_{10} + v_{6}v_{12}),
\end{eqnarray*}
where $v_i=u_i-1/2,i=1,\cdots,16$. 
\section{Discussions}\label{sec:discussions}
To conclude this paper, we address the following comments and discussions, some of which will be important for 
designs of good practical error-correcting codes.

\subsection{Stability}
Codeword fixed points in $\mathcal{C}$ are stable from Proposition \ref{pp:stable}, 
and non-codeword poles in $\mathbb{F}^n_2\setminus\mathcal{C}$ are unstable from Proposition \ref{pp:gradient} 
in the sense that nearby points to a pole in $\mathbb{F}^n_2\setminus\mathcal{C}$ leave away from the pole. 
Let us recall that each $n$-bit received sequence 
$y\in\mathbb{F}^n_2$ and its initial point $u^0\in I^n$ are related by (\ref{eq:ip}) and it is characterized that 
$y$ is the closest point to $u^0$ in $(\partial I)^n$.
Hence, if the received sequence $y$ is a codeword, then $u^0$ may approach to the codeword fixed point $y$.
This obviously depends on whether or not $u^0$ is located in the attractor region of the codeword fixed point, 
although it is actually the case when $\epsilon$ is small enough because of the stability. 
So far, a general structure of the attractor region for each codeword fixed point is not yet known.  
However this is an important subject since it is indispensable to give an estimate of error probabilities
of ML decoding and its approximation. 
Similar arguments also hold for a non-codeword received sequence and its repelling property. 
%
\subsection{Local dynamics around $p$}\label{sec:p}
Let us recall that the ML decoding rule is given by (\ref{eq:mlrule}) which checks the location of the image for 
an initial point $u^0$ to the point $p$. 
Hence the local dynamics around the point $p$ will be important for decoding process.
In the following, we explain the local dynamics around $p$ in two different cases: 
the Jacobi matrix $J$ at $p$ is (i) identity or (ii) not identity.

In the case (i), the local dynamics around $p$ is precisely determined by Theorem \ref{thm:thm}.
As explained after Theorem \ref{thm:thm}, the nonlinear dynamics around $p$ is closely related to 
the encoding structure of the code and the decoding process. 

In the case (ii), let us suppose that $B_i$ in Corollary \ref{corollary:eigenvalue} 
induces unstable eigenvalues $n_i>1$ and let us focus on its stable/unstable eigenspaces $S_i/U_i$, respectively.
From its eigenvector, $U_i$ is given by the 1 dimensional subspace spanned by
\[
(0,\cdots,0,\underbrace{1,\cdots,1}_{n_i},0,\cdots,0).
\]
This plays a role to make $u_{i_1},\cdots,u_{i_{n_i}}$ to be equal, and it reflects the fact $g_{i_1}=\cdots=g_{i_{n_i}}$. 
It means that the unstable subspace $U_i$ points to codewords' directions, 
i.e., $\mathcal{C}$.
On the other hand, $S_i$ is spanned by the stable eigenvectors $p^{(i)}_2,\cdots,p^{(i)}_{n_i}$ 
in Corollary \ref{corollary:eigenvector}. 
Contrary to the unstable eigenvector, these eigenvectors play a role to generate different elements in 
$u_{i_1},\cdots,u_{i_{n_i}}$ and point to non-codewords' directions under time reversal, 
i.e., $\mathbb{F}^n_2\setminus\mathcal{C}$. 
Therefore, the fixed point $p$ can be regarded as an indicator to codewords in the sense that 
non-codeword elements shrink and codeword elements expand around $p$. 

In both cases, further nonlinear analysis of center/stable/unstable manifolds of $p$ will be useful for 
finding suitable encoding rules and estimating the decoding performance for the approximate ML.

\subsection{Hyperbolicity of $p$ and rate restriction}
From the above argument on the fixed point $p$, it seems to be appropriate to design a generator matrix
to be hyperbolic at $p$, 
because $p$ separates expanding and shrinking directions properly and these separations have an affect 
on the decoding performance. 
However, if $p$ is hyperbolic, then the coding rate must satisfy $r=k/n\leq 1/2$ by Corollary \ref{corollary:eigenvalue}. 
Namely, the hyperbolicity prevents a code to have a high coding rate greater than half, 
although this is not a strict restriction in particular applications like wireless communication channels. 
Therefore it is necessary to have a center eigenspace at $p$ for a code with the rate $r>1/2$.
\subsection{Normal form theory in dynamical systems}
The normal form theory in dynamical systems (e.g., see \cite{wiggins}) enables us to 
transform a map into a simpler form by using a near identity transformation around a fixed point. 
One of the essential points is that nonresonant higher order terms can be removed from the original map
by this transformation. 
Theorem \ref{thm:thm} can be interpreted from the viewpoints of normal forms in such a way that
an algebraic geometry code gives only zero nonresonant terms in the expansion of $f(u)$ at $p$. 
Then, it leads to the following natural question whether a code whose rational map 
does not have resonant terms, but has nonresonant terms which are not necessarily zeros
in its expansion is a good error-correcting code or not through a near identity transformation.
At least, this class of codes contains algebraic geometry codes as a subclass, and
a similar statement to Corollary \ref{corollary:ag} holds through near identity transformations.
\subsection{Relation to LDPC codes}
It seems to be valuable to mention a relationship to LDPC codes \cite{GallagerLDPC}, \cite{MacKay},
which are a relatively new class of error-correcting codes based on iterative decoding schemes 
(for a reference to this research region, see \cite{RU}). 
The iterative decoding schemes  mainly use so-called sum-product algorithm for ML decoding and 
deal with a marginalized conditional probability in (\ref{eq:bitml}) as a convergent point.
Although this coding scheme gives a good performance in some numerical simulations, 
mathematical further understanding of the sum-product algorithm and ML decoding is 
desired to design better coding schemes. 
From the viewpoint of dynamical systems, it seems to be natural to formulate the sum-product algorithm
or ML decoding itself as a certain map, and then analyze its mechanism. 
The strategy in this paper is based on this consideration.
\section*{Acknowledgment} 
The authors express their sincere gratitude to the members of TIN working group for valuable
comments and discussions on this paper. This work is supported by JST PRESTO program.


\begin{thebibliography}{99}

\bibitem{fulton}
W. Fulton, Algebraic Curves: An Introduction to Algebraic Geometry, Addison Wesley Publishing Company, 1989.

\bibitem{Gallagerinfo}
R. G. Gallager, Information Theory and Reliable Communication, John Wiley and Sons, 1968.

\bibitem{GallagerLDPC}
R. G. Gallager, Low-Density Parity-Check Codes. Cambridge, MA: M.I.T. Press, 1963.


\bibitem{vanlint}
T. H{\o}holdt, J.H. van Lint, and R. Pellikaan, Algebraic geometry codes, V.S. Pless, W.C. Huffman (Eds.), Handbook of Coding Theory, vol. 2, Elsevier, Amsterdam, 1998, pp. 871-961.

\bibitem{MacKay}
D. J. C. MacKay, Good Error-Correcting Codes Based on Very Sparse Matrices, IEEE Trans. Inform. Theory, vol. 45, 
pp. 399-431, 1999.

\bibitem{RU}
T. Richardson and R. Urbanke, Modern Coding Theory, Cambridge University Press, 2008.

\bibitem{shannon}
C. E. Shannon, A Mathematical Theory of Communication, Bell System Technical Journal, vol. 27, pp. 379-423 and 623-656, 1948.

\bibitem{stich}
H. Stichtenoth, Algebraic Function Fields and Codes, Graduate Texts in Mathematics 2nd ed., Springer, 2008.


\bibitem{tvn}
M. Tsfasman, S. Vl\v{a}du\c{t}, and D. Nogin, Algebraic Geometric Codes: Basic Notions, Mathematical Surveys and Monographs, vol. 139, AMS.

\bibitem{wiggins}
S. Wiggins, Introduction to Applied Nonlinear Dynamical Systems and Chaos, Texts in Applied Mathematics 2, Springer-Verlag, 2003.

\end{thebibliography}
\end{document}